\RequirePackage[l2tabu,orthodox]{nag}
\pdfoutput=1

\newif\ifarxiv
\documentclass[format=acmsmall,timestamp=false,screen]{acmart}
\arxivtrue
\ifarxiv
\setkeys{acmart.cls}{nonacm}
\newcommand{\dx}{\ensuremath\,\mathrm{d}x}
\newcommand{\ds}{\ensuremath\,\mathrm{d}s}
\usepackage{tikz}
\usetikzlibrary{patterns}
\usetikzlibrary{calc,arrows,shapes,automata,backgrounds,petri}
\usepackage{pgfplots}
\usepackage{pgfplotstable}
\usepackage{mathtools}
\usepackage{subcaption}

\newcommand{\tr}{\operatorname{tr}}
\renewcommand{\div}{\operatorname{div}}

\newcommand{\bff}{\mathbf{f}}

\newcommand{\hdiv}{\ensuremath{H(\operatorname{div})}}
\newcommand{\hcurl}{\ensuremath{H(\operatorname{curl})}}
\newcommand{\HCT}{\ensuremath{\mathrm{HCT}}}
\renewcommand{\P}{\ensuremath{\mathrm{P}}}

\def \bfx{\mathbf{x}}
\def \bfv{\mathbf{v}}

\def \bfs{\mathbf{s}}
\def \bfe{\mathbf{e}}
\def \bfn{\mathbf{n}}
\def \bft{\mathbf{t}}

\usepackage{listings}
\definecolor{DarkBlue}{rgb}{0.00,0.00,0.55}
\definecolor{DarkRed}{rgb}{0.55,0.00,0.00}
\definecolor{DarkGreen}{rgb}{0.00,0.55,0.00}
\definecolor{Bittersweet}{rgb}{1.0, 0.44, 0.37}
\definecolor{Purple}{rgb}{0.5, 0.0, 0.5}

\acmJournal{TOMS}

\title{FIAT: enabling classical and modern macroelements}

\author{Pablo D. Brubeck}
\affiliation{%
  \institution{University of Oxford}
  \department{Mathematical Institute}
  \city{Oxford}
  \country{UK}}
\email{brubeckmarti@maths.ox.ac.uk}

\author{Robert C. Kirby}
\affiliation{%
  \institution{Baylor University}
  \department{Department of Mathematics}
  \streetaddress{1410 S.~4th St.}
  \city{Waco}
  \state{TX}
  \country{USA}}
\email{robert_kirby@baylor.edu}

\numberwithin{equation}{section}
\usepackage{cleveref}
\crefname{algorithm}{Algorithm}{Algorithms}
\crefname{figure}{Fig.}{Figs.}
\crefname{table}{Table}{Tables}

\begin{abstract}
  Many classical and modern finite element spaces are derived by dividing each computational cell into finer pieces.
  Such \emph{macroelements} frequently enable the enforcement of mathematically desirable properties such as divergence-free conditions or $C^1$ continuity in a simpler or more efficient manner than elements without the subdivision.
  Although a few modern software projects provide one-off support for particular macroelements, a general approach facilitating broad-based support has, until now, been lacking.
  In this work, we describe a major addition to the FIAT project to support a wide range of different macroelements.
  These enhancements have been integrated into the Firedrake code stack.
  We provide numerical evaluation of the new macroelement facility.
\end{abstract}

\begin{document}
\maketitle

\section{Introduction}
Finite element methods provide a powerful suite of tools for the numerical approximation of solutions to partial differential equations posed on $H^1$, $\hcurl$, $\hdiv$, and other Sobolev spaces.
Finite element methods have broad applicability, including in unstructured geometry, and allow general approximation orders.
At the same time, realizing this theoretical flexibility presents a technical challenge to designing general and efficient software.
FIAT, the FInite element Automatic Tabulator, was first introduced
some two decades ago to provide a general tool for the purpose of realizing this theoretical flexibility in practice \citep{Kirby:2004}.
FIAT is an independent library providing a suite of cells, basis functions, and quadrature rules, usable in principle by any client.  
After its introduction in~\citep{Kirby:2004}, it was updated to recast internal operations in terms of dense linear algebra in~\citep{kirby2006optimizing}.
FIAT's interaction with code generation for variational forms was pursued in papers such as~\citep{kirby2006compiler, rognes2010efficient}.
Support for constructing new elements via tensor products and other higher-order operations was developed in~\citep{mcrae2016automated}.
FInAT~\citep{homolya2017exposing, finat-zany} wrapped abstract syntax around FIAT, enabling code generation of more complex algorithms, including vectorized and sum-factored structured algorithms in~\citep{homolya2018tsfc} and nonstandard pullbacks required by the theory in~\citep{aznaran2022transformations, kirby-zany, bock2024planar}.

\emph{Macroelements} have long been absent from the vast array of elements supported in FIAT and most other general finite element libraries.
Based on piecewise polynomials over some subdivision of each computational cell, such elements provide important properties with lower polynomial degree and/or fewer degrees of freedom than the pure polynomial spaces previously supported in FIAT.
For example, obtaining $C^1$ continuity on triangles with pure polynomial spaces requires at least degree five.
The Argyris element uses all quintic polynomials and has 21 degrees of freedom, and the Bell element has 18 degrees of freedom and contains all quartics but not quite all quintics (the normal derivative on the edges is constrained by be cubic rather than quartic).
However, the Hsieh--Clough--Tocher (HCT) macroelement~\cite{clough1965finite} requires only 12 degrees of freedom by splitting the triangle at its barycenter and using piecewise cubics.
These elements have recently been extended to higher order variants~\citep{groselj22}.
Unlike Bell and Argyris, these elements do not employ second derivatives at vertices as degrees of freedom, which may motivate their use at higher order.

Macroelements also permit low-degree pointwise divergence-free approximations to incompressible flow.
For example, the Scott--Vogelius pair~\cite{BreSco,scott1984conforming} discretizes the Stokes equations with conforming velocities of at least quartic degree and discontinuous pressures of one degree lower.  Tetrahedra require at least hexic velocities.
While, modulo mild mesh restrictions, this gives a stable approximation with pointwise divergence-free velocities, on barycentrically refined meshes it is sufficient to use velocities with degree equal to the spatial dimension and discontinuous pressures on degree lower~\cite{guzman2018inf}.
As an alternative to modifying the given mesh, one may also implement such methods as  macroelements over a barycentric splitting.
Iso-type elements, where one uses a piecewise $\P_1$ space on a uniform refinement of each cell, also provide a use case for macroelements.
For Stokes flow, one may replace the quadratic velocities in the Taylor--Hood pair with piecewise $\P_1$ elements using the locations of the quadratic degrees of freedom, as in \Cref{p1isop2}.
Using finer refinements is known to give excellent preconditioners for high-order discretizations on quadrilateral and hexahedral meshes~\cite{pazner2023low}, although it is not known how to optimally adapt this approach to simplices.
Other novel elements, such as the Alfeld-Sorokina~\cite{alfeld-sorokina16} and Guzm\'an-Neilan~\cite{guzman2018inf} elements, can be placed in discrete complexes, and are also of current research interest.

The discretization of symmetric tensors, whether in stress-based formulations of fluids or in the Hellinger--Reissner formulation of elasticity, also motivates macroelements.
To obtain a non-macro, polynomial, $\hdiv$-conforming symmetric stress tensor requires the lowest-order Arnold--Winther triangle~\cite{Arnold2002} with 24 degrees of freedom.  On tetrahedra, the conforming Arnold-Awanou-Winther element~\cite{Arnold2008} has a hefty 162 degrees of freedom.
Somewhat smaller but nonconforming elements are known~\cite{Arnold2003,Arnold2014}.
On the other hand,~\cite{gopalakrishnan2024johnson} gives a macroelement in any dimension utilizing only piecewise linear polynomials on a barycentric refinement.  The two-dimensional element is conforming with 15 degrees of freedom, and the three-dimensional element only has 42.

Here, we extend FIAT to enable a broad class of simplicial macroelements and integrate them with the Firedrake project~\cite{rathgeber2016firedrake, FiredrakeUserManual}.
Other high-level packages support  certain basic macroelements.
For example, Deal.II and GetFem++~\cite{bangerth2007deal, renard2020getfem} support iso-type elements and the HCT triangle.  Freefem~\cite{freefem} does as well, although iso-type elements require constructing nested meshes.
Libmesh also supports HCT and some Powell-Sabin splines~\cite{stogner2007c1}.
The basix library~\cite{scroggs2022basix} has enabled iso-type elements, but these are not supported in the rest of the FEniCS code stack.
Now, FIAT goes far beyond these one-off implementations, with a much wider suite of macroelements, a more general facility for implementing new ones, and a clean integration with the rest of Firedrake.

The rest of the paper is organized as follows.  In Section~\ref{sec:macro}, we describe our general approach and a suite of macroelements enabled by the new technology.
This technology is further described in Section~\ref{sec:code}, where we detail modifications to FIAT to support macroelements and necessary changes to the remainder of the Firedrake code stack to enable their seamless use.
Most of these elements do not map simply via affine or Piola pullbacks, and
we describe the application of the theory in~\cite{kirby-zany} to transform the HCT elements in Section~\ref{sec:xform}.  Similar techniques hold for the other macroelements we have implemented.
Finally, numerical results evaluating the newly-enabled elements are given in Section~\ref{sec:num} before concluding thoughts presented in Section~\ref{sec:conc}.

\section{Macroelements}
\label{sec:macro}
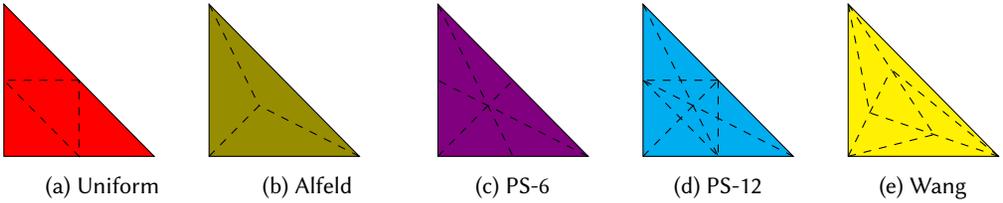
\begin{figure}[htbp]
  \centering
  \begin{subfigure}[t]{0.19\textwidth}
  \begin{tikzpicture}
    \coordinate (v0) at (0,0);
    \coordinate (v1) at (2,0);
    \coordinate (v2) at (0,2);
    \coordinate (e0) at ($0.5*(v1)+0.5*(v2)$);
    \coordinate (e1) at ($0.5*(v0)+0.5*(v2)$);
    \coordinate (e2) at ($0.5*(v0)+0.5*(v1)$);
    \draw[fill=red] (v0) -- (v1) -- (v2) -- cycle;
    \draw[dashed] (e0) -- (e1) -- (e2) -- cycle;
  \end{tikzpicture}
  \caption{Uniform}
  \label{fig:unif}
  \end{subfigure}
  \centering
  \begin{subfigure}[t]{0.19\textwidth}
  \begin{tikzpicture}
    \coordinate (v0) at (0,0);
    \coordinate (v1) at (2,0);
    \coordinate (v2) at (0,2);
    \coordinate (b) at ($0.333*(v0)+0.333*(v1)+0.333*(v2)$);
    \draw[fill=olive] (v0) -- (v1) -- (v2) -- cycle;
    \draw[dashed] (v0) -- (b);
    \draw[dashed] (v1) -- (b);
    \draw[dashed] (v2) -- (b);
  \end{tikzpicture}
  \caption{Alfeld}
  \label{fig:alfeldsplit}
  \end{subfigure}
  \begin{subfigure}[t]{0.19\textwidth}
  \centering
  \begin{tikzpicture}
    \coordinate (v0) at (0,0);
    \coordinate (v1) at (2,0);
    \coordinate (v2) at (0,2);
    \coordinate (b) at ($0.333*(v0)+0.333*(v1)+0.333*(v2)$);
    \coordinate (e0) at ($0.5*(v1)+0.5*(v2)$);
    \coordinate (e1) at ($0.5*(v0)+0.5*(v2)$);
    \coordinate (e2) at ($0.5*(v0)+0.5*(v1)$);
    \draw[fill=violet] (v0) -- (v1) -- (v2) -- cycle;
    \draw[dashed] (v0) -- (e0);
    \draw[dashed] (v1) -- (e1);
    \draw[dashed] (v2) -- (e2);
  \end{tikzpicture}
  \caption{PS6}
  \label{fig:pssplit}
  \end{subfigure}
  \begin{subfigure}[t]{0.19\textwidth}
  \centering
  \begin{tikzpicture}
    \coordinate (v0) at (0,0);
    \coordinate (v1) at (2,0);
    \coordinate (v2) at (0,2);
    \coordinate (b) at ($0.333*(v0)+0.333*(v1)+0.333*(v2)$);
    \coordinate (e0) at ($0.5*(v1)+0.5*(v2)$);
    \coordinate (e1) at ($0.5*(v0)+0.5*(v2)$);
    \coordinate (e2) at ($0.5*(v0)+0.5*(v1)$);
    \draw[fill=cyan] (v0) -- (v1) -- (v2) -- cycle;
    \draw[dashed] (v0) -- (e0);
    \draw[dashed] (v1) -- (e1);
    \draw[dashed] (v2) -- (e2);
    \draw[dashed] (e0) -- (e1);
    \draw[dashed] (e1) -- (e2);
    \draw[dashed] (e0) -- (e2);
  \end{tikzpicture}
  \caption{PS12}
  \label{fig:pssplit12}
  \end{subfigure}
  \begin{subfigure}[t]{0.19\textwidth}
    \centering
    \begin{tikzpicture}[scale=0.667]
      \draw[fill=yellow] (0,0) -- (3,-3) -- (0,-3) -- cycle;
      \coordinate (v0) at (0, 0);
      \coordinate (v1) at (3, -3);
      \coordinate (v2) at (0, -3);
    \coordinate (w0) at ($0.571*(v0)+0.286*(v1)+0.143*(v2)$);
    \coordinate (w1) at ($0.571*(v1)+0.286*(v2)+0.143*(v0)$);
    \coordinate (w2) at ($0.571*(v2)+0.286*(v0)+0.143*(v1)$);
    \draw[dashed] (v0) -- (w1);
    \draw[dashed] (v0) -- (w2);
    \draw[dashed] (v1) -- (w0);
    \draw[dashed] (v1) -- (w2);
    \draw[dashed] (v2) -- (w0);
    \draw[dashed] (v2) -- (w1);
    \end{tikzpicture}
    \caption{Wang}
    \label{fig:wang}
  \end{subfigure}
  \caption{Some typical splitting strategies for macroelements.}
  \Description{Pictures of split triangles.}
  \label{fig:splits}
\end{figure}

Typical finite element spaces employ functions that are piecewise polynomial over each cell in the computational domain, with some restrictions on the continuity between cells.
However, macroelements require each cell to be further subdivided in some regular way, so that the functions on each cell in the mesh are themselves piecewise polynomials.
\Cref{fig:splits} depicts several common splittings used on triangles.
Many of these splitting strategies, such as the uniform split in \Cref{fig:unif} and Alfeld (also known as barycentric or Clough--Tocher splitting) in \Cref{fig:alfeldsplit} are affinely invariant -- affinely mapping a split reference cell and yields the correctly split physical cell.
The Powell-Sabin-12 split in \Cref{fig:pssplit12} and Wang splitting in \Cref{fig:wang} may also be constructed in this way.

Certain splittings, such as the 6-way Powell-Sabin split in \Cref{fig:pssplit}, present some subtlety.
Following the argument in~\cite{powell1977piecewise}, one selects a point in the interior of each triangle.  The line segment connecting interior points of adjacent (i.e. sharing an edge) triangles must intersect the interior of the triangles' shared edge in order to guarantee the existence of $C^1$ piecewise quadratic splines.
The line connecting incenters of adjacent triangles always possess this property, so that the construction works on general meshes.
This can only fail if the union of adjacent triangles is nonconvex (and even then, it is not automatic), so barycentric splitting should be sufficient on quality triangulations.
Similar considerations arise for the Worsey--Farin and Worsey--Piper splits of tetrahedra~\cite{lai2007spline, worsey1988trivariate}.
Constructions on general meshes also require geometrically-dependent points that in general, quite difficult to compute.  Whether mesh restrictions imposed by barycentric splitting would be as gentle as for triangles, however, remains unclear.
At any rate, the use of a reference element is deeply ingrained in general purpose finite element codes such as Firedrake.  Our development concentrates on affinely-preserved splits based on barycenters, which seems perfectly acceptable in two dimensions at least.

Before describing our implementation of macroelements in FIAT and the rest of the Firedrake code stack, we give a careful description of the various elements depicted above.  In line with the rest of FIAT, we work in terms of the Ciarlet triple~\cite{ciarlet2002finite}.
\begin{itemize}
\item $K \subset \mathbb{R}^d$ is a bounded domain with piecewise smooth boundary.  In this paper, we only consider $K$ as a simplex.
\item $P$ is a finite-dimensional function space defined on the closure of $K$, typically consisting of polynomials or vectors/tensors of them, or piecewise polynomials over a subdivision of $K$.
\item $N = {\left\{ n_i \right\}}_{i=1}^{\dim{P}}$ is a basis for the dual space $P^\prime$, called the set of \emph{nodes} or \emph{degrees of freedom}.
\end{itemize}
The \emph{nodal basis} for a finite element is the set ${\left\{ \phi_i \right\}}_{i=1}^{\dim P} \subset P$ such that,
\begin{equation}
n_i(\phi_j) = \delta_{ij}, \ \ \ 1 \leq i, j \leq \dim P.
\end{equation}
The nodes of a finite element typically consist of functionals such as pointwise evaluation of functions or derivatives at particular points, or certain integral moments of functions on $K$ or its boundary facets, and are chosen to enforce certain kinds of continuity between adjacent elements.

In FIAT, we compute the nodal basis numerically by means of solving a generalized Vandermonde system.
Given any readily computable basis $\{ p_i \}_{i=1}^{\dim P}$ for $P$, we can write
\begin{equation}
\phi_j = \sum_{k=1}^{\dim P} A_{jk} p_k,
\end{equation}
and applying any node $n_i$ to both sides of the equation gives
\begin{equation}
\delta_{ij} = \sum_{k=1}^{\dim P} A_{jk} n_i(p_k),
\end{equation}
so that, with $V_{ij} = n_i(p_j)$,
\begin{equation}
   I = A V^\top.
\end{equation}

Putting macroelements into this framework represents a departure from the Bernstein--Bezier techniques more commonly used for representing triangular splines.
These techniques directly define the macroelement basis on each cell in terms of geometrically-dependent linear combinations of the Bernstein polynomials.
In some situations, this can allow nonnegative basis functions~\cite{groselj22}, or avoid mesh restrictions~\cite{powell1977piecewise}.
A local Bernstein representation also could 
allow sum-factored algorithms~\cite{ainsworth2011bernstein,kirby2011fast,kirby2014low}.
However, using a reference element and Vandermonde matrix
retains the declarative approach provided by FIAT -- one
simply implements an element by defining a space of functions and a list of
functionals defining the dual basis, whence the nodal basis is
constructed automatically.
While this limits our geometric flexibility, it greatly simplifies the
path to implementation, both of specific macroelements and their
incorporation in the rest of the Firedrake code stack.

\subsection{Splittings}
We let $K$ be any simplex in $\mathbb{R}^d$.
We let $\Delta_A(K)$ be the barycentric/Alfeld split of $K$ into $d+1$ simplices, shown for triangles in \Cref{fig:alfeldsplit}.
By dividing each edge into $\ell$ uniform segments and connecting, we split a triangle into $\ell^2$ subcells.  This \emph{iso-split}, denoted by
$\Delta_{iso,\ell}(K)$ is shown with $\ell=2$ in \Cref{fig:unif}.
We omit the subscript when $\ell=2$, so $\Delta_{iso}(K) \equiv \Delta_{iso,2}(K)$.
For triangles, we let $\Delta_{PS6}(K)$ and $\Delta_{PS12}(K)$ be the sets of triangles obtained via the Powell-Sabin splits of \Cref{fig:pssplit,fig:pssplit12}.

%% Each of these splittings defines a simplicial complex in $\mathbb{R}^d$.  For $0 \leq i \leq d$, we let $\Delta^i$ denote the set of entities of dimension $i$ (e.g. vertices when $i=0$).
%% Similarly, the simplex $K$ itself is a complex, and we let $K^i$ denote its entities of dimension $i$.  We will frequent

\subsection{Degrees of freedom}
In order to define macroelements, we need to establish some notation for degrees of freedom.
For each vertex $\mathbf{v} \in K$ of the unsplit simplex, we define the functional $\delta_\mathbf{v}$ by
\begin{equation}
  \delta_{\mathbf{v}}(p) = p(\mathbf{v}).
\end{equation}
Similarly, we differentiate at $\mathbf{v}$ in the direction of some unit vector $\mathbf{s}$ by
\begin{equation}
  \delta_{\mathbf{v}}^{\mathbf{s}}(p) = \mathbf{s}\cdot \nabla p(\mathbf{v}).
\end{equation}
The gradient at a point is a vector of functionals.  It takes the form
\begin{equation}
 \nabla_{\mathbf{v}} = \begin{bmatrix} \delta_{\mathbf{v}}^{\mathbf{x}} & \delta_{\mathbf{v}}^{\mathbf{y}} \end{bmatrix}^\top
\end{equation}
in two dimensions, with similar form in three.

We also need to define functionals based on integral moments over edges or faces.
For any facet $\mathbf{f}$ of $K$ itself or its splitting and any function $p \in L^2(\mathbf{f})$, we define the
integral moment
\begin{equation}
  \mu_{\mathbf{f}, q}(p) = \int_{\mathbf{f}} q p \ds.
\end{equation}
and for unit vector $\mathbf{s}$ and facet $\mathbf{f}$, we define the integral moment of the directional derivative by
\begin{equation}
\label{eq:edgederivmoment}
  \mu^{\mathbf{s}}_{\mathbf{f}, q}(p) = \int_{\mathbf{f}} q \left( \mathbf{s}\cdot \nabla p \right) \ds.
\end{equation}
With these moment-based functionals, we omit the subscript $q$ when $q \equiv 1$ on $\mathbf{f}$. For instance,
$\mu^{\mathbf{s}}_{\mathbf{f}} \equiv \mu^{\mathbf{s}}_{\mathbf{f},1}(p)$.

\subsection{$C^0$ macroelements}
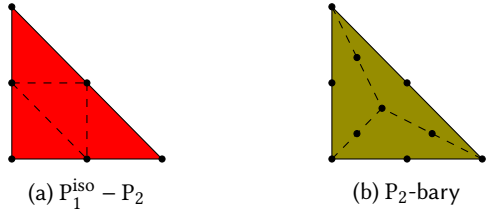
\begin{figure}[htbp]
  \begin{subfigure}[t]{0.3\textwidth}
    \centering
    \begin{tikzpicture}[scale=1.0] % p1isop2
      \draw[fill=red] (0,0) -- (2, 0) -- (0, 2) -- cycle;
      \foreach \n in {0,...,2} {
        \foreach \i in {0,...,\n} {
          \draw[fill] (\i,2-\n) circle (0.04);
        }
      }
      \draw[dashed] (0,1) -- (1,1);
      \draw[dashed] (1,1) -- (1,0);
      \draw[dashed] (0,1) -- (1,0);
    \end{tikzpicture}
    \caption{$\P_1^\mathrm{iso}-\P_2$}
    \label{p1isop2}
  \end{subfigure}
  \begin{subfigure}[t]{0.3\textwidth}
    \centering
    \begin{tikzpicture}[scale=0.667] % SV
    \draw[fill=olive] (0,0) -- (3, 0) -- (0, 3) -- cycle;
    \draw[dashed] (0,0) -- (1,1);
    \draw[dashed] (3,0) -- (1,1);
    \draw[dashed] (0,3) -- (1,1);
    \foreach \x/\y in {0/0,3/0,0/3,1/1,0/1.5,1.5/1.5,1.5/0,0.5/0.5,2/0.5,0.5/2} {
      \draw[fill] (\x, \y) circle (0.06);
      }
    \end{tikzpicture}
    \caption{$\P_2$-bary}
    \label{p2bary}
  \end{subfigure}
  \caption{Some Lagrange-type $C^0$ macroelements. Solid dots represent point evaluation degrees of freedom.}
  \Description{Schematic pictures showing degrees of freedom for Lagrange macroelements.}
  \label{fig:macrolag}
\end{figure}

For any simplex $K$, we define the space of $C^0$ splines of degree $k$ over a splitting $\Delta$ of $K$ 
\begin{equation}
  S^0_k(\Delta) = \{ s \in C^0(K) : s|_\tau \in \P_k(\tau), \tau \in \Delta \}.
\end{equation}
Then, we parametrize the space by standard degrees of freedom. If $K$ is a triangle, we choose degrees of freedom as
\begin{itemize}
\item $\delta_{\mathbf{v}}$ for each vertex $\mathbf{v}$ of $\Delta$.
\item If $k > 1$, take moments $\mu_{\bfe, q}$ for each edge $\mathbf{e}$ in $\Delta$ and each $q$ in a basis for $\P_{k-2}(\mathbf{e})$.
\item If $k > 2$, take moments $\mu_{\tau, q}$ for each $\tau \in \Delta$ and each $q$ in a basis for $\P_{k-3}(\tau)$.
\end{itemize}
Alternatively (and far more commonly for low-order finite elements), one may replace the integral moments with evaluation at an appropriate set of unisolvent points.
\Cref{fig:macrolag} shows the degree-of-freedom diagrams for some $C^0$ macroelements.

\subsection{$C^1$ macroelements}
For $C^1$ elements, we focus on splittings of a triangle $K$ rather than general simplex.  For some splitting $\Delta$ of $K$, we define
\begin{equation}
S_k^1(\Delta) = \left\{ s \in C^1(K) : s|_\tau \in \P_k(\tau), \tau \in \Delta \right\}
\end{equation}
to be the space of continuously differentiable piecewise polynomials of degree $k$.
When $k=2$ and $\Delta \in \{ \Delta_{PS6}(T), \Delta_{PS12}(T) \}$, we obtain the spaces for the quadratic Powell-Sabin splines in \Cref{fig:pssplit,fig:pssplit12}.
When $k=3$ and $\Delta = \Delta_{A}(K)$, this is the standard HCT space.

Letting $\mathcal{E}$ denote the edges of $K$, we also obtain the reduced space
\begin{equation}
\tilde{S}_3(\Delta) = \left\{ s \in S_3^1(\Delta) : \tfrac{\partial s}{\partial n}|_e \in \P_1(e), e \in \mathcal{E} \right\}.
\end{equation}
Much like the non-macro Bell triangle, this smaller space removes edge degrees of freedom by restricting the degree of the normal derivative of members of the function space on each edge.

Recently, Groselj and Knez~\cite{groselj22} have constructed a generalization of the space $S_k^1(\Delta)$ to higher-order polynomials while retaining $C^1$ continuity and full approximating power.  They define the supersmooth space
\begin{equation}
   S_k(\Delta_{A}(K)) = S^1_k(\Delta_{A}(K)) \cap C^{k-1}(\mathbf{v}_0),
\end{equation}
where $\mathbf{v}_0$ denotes the interior point of the Alfeld split $\Delta_{A}(K)$.
This space has greater than $C^1$ continuity within the triangle $K$, although it only joins across triangles with $C^1$ continuity.

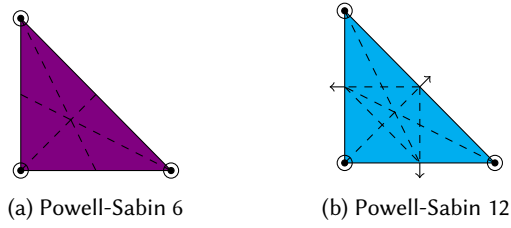
\begin{figure}[htbp]
  \begin{subfigure}[t]{0.3\textwidth}
    \centering
    \begin{tikzpicture}[scale=0.667] %PS2
    \draw[fill=violet] (0,0) -- (3, 0) -- (0, 3) -- cycle;
    \draw[dashed] (0,0) -- (1,1);
    \draw[dashed] (3,0) -- (1,1);
    \draw[dashed] (0,3) -- (1,1);
    \draw[dashed] (1.5,0) -- (1,1);
    \draw[dashed] (1.5,1.5) -- (1,1);
    \draw[dashed] (0,1.5) -- (1,1);
    \foreach \i/\j in {0/0, 3/0, 0/3}{
      \draw[fill=black] (\i, \j) circle (0.06);
      \draw (\i, \j) circle (0.15);
    }
    \end{tikzpicture}
    \caption{Powell-Sabin 6}
    \label{PS6}
  \end{subfigure}
  \begin{subfigure}[t]{0.3\textwidth}
  \centering
  \begin{tikzpicture}
    \coordinate (v0) at (0,0);
    \coordinate (v1) at (2,0);
    \coordinate (v2) at (0,2);
    \coordinate (b) at ($0.333*(v0)+0.333*(v1)+0.333*(v2)$);
    \coordinate (e0) at ($0.5*(v1)+0.5*(v2)$);
    \coordinate (e1) at ($0.5*(v0)+0.5*(v2)$);
    \coordinate (e2) at ($0.5*(v0)+0.5*(v1)$);
    \draw[fill=cyan] (v0) -- (v1) -- (v2) -- cycle;
    \draw[thin,dashed] (v0) -- (e0);
    \draw[thin,dashed] (v1) -- (e1);
    \draw[thin,dashed] (v2) -- (e2);
    \draw[thin,dashed] (e0) -- (e1);
    \draw[thin,dashed] (e1) -- (e2);
    \draw[thin,dashed] (e0) -- (e2);
    \foreach \i/\j in {0/0, 2/0, 0/2}{
      \draw[fill=black] (\i, \j) circle (0.04);
      \draw (\i, \j) circle (0.1);
    }
    \foreach \i/\j/\n/\t in {1/0.0/0.0/-2, 1/1/1.41/1.41, 0.0/1/-2/0}{
      \draw[->] (\i, \j) -- (\i+\n/10, \j+\t/10);
    }
  \end{tikzpicture}
  \caption{Powell-Sabin 12}
  \label{PS12}
  \end{subfigure}
  \caption{Quadratic $C^1$ macroelements on Powell-Sabin splits. Hollow circles represent derivative evaluation for each Cartesian direction, and
    arrows represent normal derivative moments along edges.}
  \Description{Schematic for quadratic Powell-Sabin macroelements.}
  \label{fig:psc1}
\end{figure}

Now, the nine degrees of freedom for  $S^1_2(\Delta_{PS6}(K))$, shown 
\Cref{PS6}, can be taken to be:
\begin{itemize}
\item $\delta_{\mathbf{v}}$ for each vertex $\mathbf{v}$ of $K$,
\item $\delta_{\mathbf{v}}^{\mathbf{s}}$ for each vertex $\mathbf{v}$ of $K$ and for $\mathbf{s}$ each Cartesian direction $\mathbf{x}$ and $\mathbf{y}$.
\end{itemize}

Similarly, the twelve degrees of freedom parameterizing  $S^1_2(\Delta_{PS12}(K)) $, shown in \Cref{PS12}, are:
\begin{itemize}
\item $\delta_{\mathbf{v}}$ for each vertex $\mathbf{v}$ of $K$,
\item $\delta_{\mathbf{v}}^{\mathbf{s}}$ for each vertex $\mathbf{v}$ of $K$ and for $\mathbf{s}$ each Cartesian direction $\mathbf{x}$ and $\mathbf{y}$.
\item $\mu_\bfe^\mathbf{n_e}$ for each edge $\mathbf{e}$ of $K$, where $\mathbf{n_e}$ is the normal to each edge.
\end{itemize}

\begin{figure}[htbp]
  \centering
  \begin{subfigure}[t]{0.3\textwidth}
    \centering
    \begin{tikzpicture}[scale=2.0] % hct
    \draw[fill=cyan] (0,0) -- (1, 0) -- (0, 1) -- cycle;
    \foreach \i/\j in {0/0, 1/0, 0/1}{
      \draw[fill=black] (\i, \j) circle (0.02);
      \draw (\i, \j) circle (0.05);
    }
    \draw[dashed] (0,0) -- (1/3,1/3);
    \draw[dashed] (1,0) -- (1/3,1/3);
    \draw[dashed] (0,1) -- (1/3,1/3);
  \end{tikzpicture}
     \caption{Reduced $\HCT$}
  \label{rhct}
  \end{subfigure}
  \begin{subfigure}[t]{0.3\textwidth}
    \centering
    \begin{tikzpicture}[scale=2.0] % hct
    \draw[fill=yellow] (0,0) -- (1, 0) -- (0, 1) -- cycle;
    \foreach \i/\j in {0/0, 1/0, 0/1}{
      \draw[fill=black] (\i, \j) circle (0.02);
      \draw (\i, \j) circle (0.05);
    }
    \foreach \i/\j/\n/\t in {0.5/0.0/0.0/-1, 0.5/0.5/0.707/0.707, 0.0/0.5/-1/0}{
      \draw[->] (\i, \j) -- (\i+\n/10, \j+\t/10);
    }
    \draw[dashed] (0,0) -- (1/3,1/3);
    \draw[dashed] (1,0) -- (1/3,1/3);
    \draw[dashed] (0,1) -- (1/3,1/3);
  \end{tikzpicture}
  \caption{$\mathrm{HCT}_3$}
  \label{hct}
  \end{subfigure}
  \begin{subfigure}[t]{0.3\textwidth}
    \centering
    \begin{tikzpicture}[scale=2.0] % hct
    \draw[fill=gray] (0,0) -- (1, 0) -- (0, 1) -- cycle;
    \foreach \i/\j in {0/0, 1/0, 0/1}{
      \draw[fill=black] (\i, \j) circle (0.02);
      \draw (\i, \j) circle (0.05);
    }
    \foreach \i/\j/\n/\t in {0.33/0.0/0.0/-1, 0.66/0.0/0.0/-1,
      0.33/0.66/0.707/0.707, 0.66/0.33/0.707/0.707,
      0.0/0.33/-1/0, 0.0/0.66/-1/0
    }{
      \draw[->] (\i, \j) -- (\i+\n/10, \j+\t/10);
    }
    \foreach \i/\j/\n/\t in {0.5/0.0/0.0/-1, 0.5/0.5/0.707/0.707, 0.0/0.5/-1/0}{
      \draw[very thick,-] (\i+\t/20, \j-\n/20) -- (\i-\t/20, \j+\n/20);
    }
    \draw[dashed] (0,0) -- (1/3,1/3);
    \draw[dashed] (1,0) -- (1/3,1/3);
    \draw[dashed] (0,1) -- (1/3,1/3);
  \end{tikzpicture}
  \caption{$\HCT_4$}
  \label{hct4}
  \end{subfigure}  
  \caption{HCT-type $C^1$ macroelements on the Alfeld split. Solid lines represent moments along edges.}
  \Description{Schematic for HCT-type macroelements.}  
  \label{fig:alfeldc1}
\end{figure}
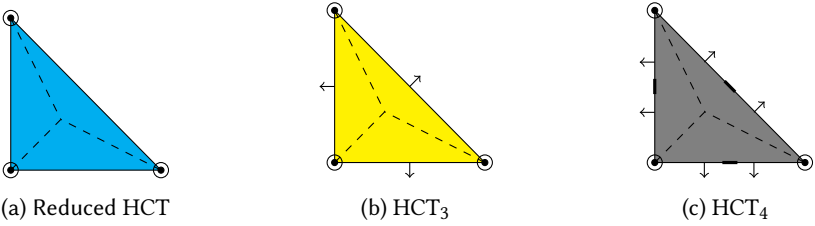

The 12-dimensional HCT space, $S^1_3(\Delta_{A}(K))$, can be parametrized in the same way as $S^1_2(\Delta_{PS12}(K))$,
and the 9-dimensional reduced HCT space $\tilde{S}_3(\Delta_{A}(K))$ 
as $S^1_2(\Delta_{PS6}(K))$.

The higher-order HCT spaces require additional edge and interior degrees of freedom.  For $S_k(\Delta_A(K))$, we take degrees of freedom
\begin{itemize}
\item $\delta_{\mathbf{v}}$ for each vertex $\mathbf{v}$ of $K$,
\item $\delta_{\mathbf{v}}^{\mathbf{s}}$ for each vertex $\mathbf{v}$ of $K$ and for $\mathbf{s}$ each Cartesian direction $\mathbf{x}$ and $\mathbf{y}$.
\item $\mu_{\bfe,q^1}^{\mathbf{n}}$ for each edge $\mathbf{e}$ of $K$ and $q^1$ in a basis for $\P_{k-3}(\mathbf{e})$.  Here, $\mathbf{n}$ is the normal to edge $\mathbf{e}$
\item If $k > 3$, take the additional moments:
  \begin{itemize}
  \item $\mu_{\bfe,q^2}$, for each edge $\mathbf{e}$ of $K$ and $q^2$ in a basis for $\P_{k-4}(\mathbf{e})$.
  \item $\mu_{K, q^3}$ for each $q^3$ in a basis for $\P_{k-4}(K)$.
  \end{itemize}
\end{itemize}

We have made specific choices of polynomial bases for edge and interior moments
to expedite the mapping from a reference element and also, incidentally, give a hierarchical basis.
The normal derivative moments and moments of the trace are taken against complementary functions. In
particular, $q_2$ should be taken as the derivatives of the polynomials in $q_1$.
We take the polynomials $q^1$ on the edge $\mathbf{e}$ to be the Jacobi polynomials
$P_i^{(1,1)}$ for $0 \leq i \leq k-3$ mapped to that edge.  Then, $q^2$ are the
derivatives of the Jacobi polynomials $\tfrac{\mathrm{d}}{\mathrm{d}s} P^{(1,1)}_{i}(s)$ for $1
\leq i \leq k-3$.
Finally, we take the polynomials $q^3$ on a triangle to be the Dubiner
polynomials of degree $k-4$, ordered hierarchically so that all polynomials of
one degree occur before any of the next higher degree.

\subsection{Stokes elements}
The Stokes equations, widely studied from physical and numerical perspectives, model creeping flow of an incompressible viscous fluid.
We consider both velocity/pressure and stress/velocity formulations, described below in Equations~\eqref{eq:velprestokes} and~\eqref{eq:stressvelstokes}.
Both formulations solve for a pair of unknowns, which must be taken
from a pair of spaces satisfying an inf-sup condition.
Moreover, finding pairs of spaces for which the discrete velocity is
divergence free \emph{pointwise} is an additional challenge.
As cited above, the seemingly obvious continuous velocity and discontinuous pressures require quite high order polynomials, greatly reduced 
on barycentric refinement of an existing simplicial mesh, one needs only $k\geq d$~\cite{guzman2018inf}.
This significantly reduces the total number of degrees of freedom.  In two dimensions, each quartic velocity component requires 15 degrees of freedom
per triangle, while quadratics on the split mesh require only 10, as shown in \Cref{p2bary}.

%% TODO: is this needed here? 
%% Implementing lower order Scott--Vogelius via barycentric refinement of a given mesh is certainly possible, but presents certain drawbacks.
%% Beyond requiring additional user intervention, one must discretize any additional unknowns in a coupled problem on the refined mesh or else implement more complex logic to evaluate operators using both the original and refined mesh.
%% Instead, we implement the Scott--Vogelius pair as macro-element on the original mesh, using the \lstinline{'alfeld'} variant for the Lagrange elements.
%% No other user-level intervention is required.

Other pairs with macroelements can give pointwise divergence-free velocities with continuous pressures, which allows smaller global approximating spaces.
These lie in a differential complex and also ensure the
divergence-free condition holds pointwise, but this work seems to be their first
realization in a general-purpose software package.  In particular, we consider the Guzm\'an--Neilan and Alfeld--Sorokina macroelements. 

It is possible to obtain a minimal inf-sup stable Stokes pair that includes
the linears on the unsplit cell, plus some face bubbles.
The motivation behind enriching with face bubbles is to control the divergence
by adding normal degrees of freedom on each face.
A non-macroelement example is the Bernardi--Raugel element, which is defined
by enriching linear vector fields with normal face bubbles \cite{bernardi-raugel1985}.
It achieves inf-sup stability when paired with piecewise constants.
Nevertheless, the face bubbles are of degree $d$, and their divergence is of degree $d-1$, so
the divergence-free constraint can only be weakly enforced when testing against piecewise constants.  That is, the integral of the divergence over each cell vanishes in this case.

\begin{figure}[htbp]
  \centering
  \begin{subfigure}[t]{0.3\textwidth}
    \centering
    \begin{tikzpicture}[scale=2.0] % guzman neilan
    \draw[fill=red] (0,0) -- (1, 0) -- (0, 1) -- cycle;
    \foreach \i/\j in {0/0, 1/0, 0/1}{
      \draw[->] (\i, \j) -- (\i+0.1, \j);
      \draw[->] (\i, \j) -- (\i, \j+0.1);
    }
    \foreach \i/\j/\n/\t in {0.5/0.0/0.0/-1, 0.5/0.5/0.707/0.707, 0.0/0.5/-1/0}{
      \draw[->] (\i, \j) -- (\i+\n/10, \j+\t/10);
    }
    \draw[dashed] (0,0) -- (1/3,1/3);
    \draw[dashed] (1,0) -- (1/3,1/3);
    \draw[dashed] (0,1) -- (1/3,1/3);
  \end{tikzpicture}
  \caption{Guzm\'an--Neilan $H^1$}
  \label{gn}
  \end{subfigure}
  \begin{subfigure}[t]{0.3\textwidth}
    \centering
    \begin{tikzpicture}[scale=2.0] % alfeld sorokina
    \draw[fill=green] (0,0) -- (1, 0) -- (0, 1) -- cycle;
    \foreach \i/\j/\n/\t in {0/0/-0.707/-0.707, 1/0/0/-0.707, 0/1/-1.414/0}{
       \node[align=center] at (\i+\n/10, \j+\t/10) {\small div};
    }
    \foreach \i/\j in {0/0, 1/0, 0/1, 0.5/0, 0.5/0.5, 0/0.5}{
      \draw[->] (\i, \j) -- (\i+0.1, \j);
      \draw[->] (\i, \j) -- (\i, \j+0.1);
    }
    \draw[dashed] (0,0) -- (1/3,1/3);
    \draw[dashed] (1,0) -- (1/3,1/3);
    \draw[dashed] (0,1) -- (1/3,1/3);
  \end{tikzpicture}
     \caption{Alfeld--Sorokina $H^1(\div)$}
  \label{alfeldsorokina}
  \end{subfigure}
  \begin{subfigure}[t]{0.3\textwidth}
    \centering
    \begin{tikzpicture}[scale=2.0] % johnson mercier
    \draw[fill=magenta] (0,0) -- (1, 0) -- (0, 1) -- cycle;
    \foreach \i/\j/\n/\t in {0.33/0.0/0.0/-1, 0.66/0.0/0.0/-1,
      0.33/0.66/0.707/0.707, 0.66/0.33/0.707/0.707,
      0.0/0.33/-1/0, 0.0/0.66/-1/0
    }{
      \draw[thick,->] (\i, \j) -- (\i+\n/10, \j+\t/10);
    }
    \foreach \i/\j in {0.44/0.11, 0.11/0.44, 0.44/0.44}{
      \draw[fill=black] (\i, \j) circle (0.02);
      }
    \draw[dashed] (0,0) -- (1/3,1/3);
    \draw[dashed] (1,0) -- (1/3,1/3);
    \draw[dashed] (0,1) -- (1/3,1/3);
  \end{tikzpicture}  
  \caption{Johnson--Mercier $H(\div, \mathbb{S})$}
  \label{johnsonmercier}
  \end{subfigure}
  \caption{Lowest-order macroelements for the Stokes and elasticity complexes. 
   The thin arrows in Guzm\'an--Neilan and Alfeld--Sorokina represent evaluation of vector components.
   The thicker arrows in Johnson--Mercier represent moments of the normal-normal and normal-tangential components of a tensor.}
  \Description{Schematic for Stokes and elasticity macroelements.}
\end{figure}
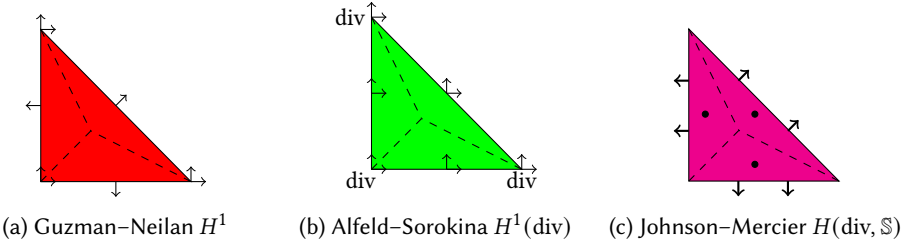

In two dimensions, Arnold and Qin \cite{arnold-qin1992} proposed a macroelement
on the Alfeld split with divergence-free quadratic face bubbles.  This
element pairs with piecewise constants on the unsplit mesh, and may be regarded
as a modification of Bernardi--Raugel. Here, the face bubbles are modified by
subtracting a piecewise polynomial on the split cell that also vanishes on the
boundary, matching the divergence of the Bernardi--Raugel face bubble.  More
recently, this construction has been extended to any dimension by Guzm\'an and
Neilan \cite{guzman2018inf}. \Cref{gn} shows the degrees of freedom of
this element on the triangle.

Macroelements also enable divergence-free formulations for Stokes with 
continuous pressure elements.  These are naturally posed in $H^1(\div) \times
H^1$, where $H^1(\div) = \{\mathbf{v} \in H^1 : \div\mathbf{v} \in H^1\}$.  The
quadratic Alfeld--Sorokina macroelement \cite{alfeld-sorokina16} is constructed
as the piecewise quadratics on an Alfeld split with $C^0$ divergence.  This
construction works in any dimension, but is only inf-sup  stable in 2D.
It is also shown in \cite{guzman2018inf} that an inf-sup stable Stokes element
with $C^0$ divergence can be obtained in three dimensions by enriching the
Alfeld-Sorokina quadratic macroelement with the Guzm\'an--Neilan cubic divergence-free face bubbles.
\Cref{alfeldsorokina} shows the degrees of freedom for this
element in two dimensions, which are those of unsplit quadratic vectors plus
the divergence at each vertex.

\subsection{The Johnson--Mercier element}
The elements of Johnson and Mercier~\cite{johnson1978some} provide symmetric $\hdiv$-conforming tensors suitable for the discretization of the Hellinger--Reissner formulation of elasticity or stress-velocity formulations of incompressible flow.
Recently, a general formulation for simplicial elements in $\mathbb{R}^d$ has been given for all $d \geq 2$~\cite{gopalakrishnan2024johnson}.
Letting $\mathbb{S}^d$ denote the space of symmetric $d \times d$ tensors, they discretize $H(\operatorname{div}; \mathbb{S}^d)$.

Let $K$ be a simplex in $\mathbb{R}^d$, split into $d+1$ subcells via the Alfeld split.
We let $\Sigma_h(K)$ denote the div-conforming space of all functions mapping $K$ into $\mathbb{S}^d$:
\[
\Sigma_h(K) = \{ \tau \in H(\operatorname{div}, K, \mathbb{S}^d):
\tau|_{K_i} \in \P_1(K_i, \mathbb{S}^d), 1 \leq i \leq d+1 \},
\]
where $\P_1(K_i, \mathbb{S}^d)$ just comprises symmetric tensor-valued functions whose components are all linear polynomials over subcell $K_i$.

This space has dimension $(d+\tfrac{1}{2})d(d+1)$, which is 15 for triangles
and 42 for tetrahedra and hence much smaller than the Arnold-Winther elements.
\Cref{johnsonmercier} shows the Johnson--Mercier element on a triangle.
Degrees of freedom parametrizing space come in two kinds:
\begin{itemize}
\item For each facet of codimension 1, $\bff$, with unit normal $\bfn$, the integral moments of each component of $\tau \bfn$ against all linear polynomials over $\bff$.
\item The integral average of each independent component of $\tau$ over $K$.
\end{itemize}

\subsection{Abbreviations}
We have introduced many finite elements, and it will be helpful, especially in tables and figures, to refer to them in shorthand.
For reference, we expand these abbreviations for Stokes pairs in Table~\ref{table:abbrev}.

\begin{table}
  \begin{tabular}{c|c|c}
    Abbreviation & Element pair & Formulation \\ \hline
    TH & Taylor--Hood & velocity/pressure \\
    ISO & $\P_2^\mathrm{iso}-\P_1$ & velocity/pressure \\
    SV & $\P_2-\P_1^{disc}$ on Alfeld split & velocity/pressure \\
    JM & Johnson--Mercier - $\P_1^{disc}$ & stress/velocity \\
    GN & Guzm\'an--Neilan - $\P_0^{disc}$ & velocity/pressure \\
    AS & Alfeld--Sorokina - $\P_1$ on Alfeld split & velocity-pressure \\ \hline
  \end{tabular}
  \caption{Abbreviations for Stokes pairs used in convergence tests}
  \label{table:abbrev}
\end{table}

\subsection{Computational costs of macroelements}
Despite providing particular features with fewer global degrees of freedom, computing with macroelements incurs some additional expense in element-level computation.
 Here, we enumerate what the costs are for a standard-fare algorithm that will iterate over all pairs of basis functions and over all quadrature points to form some local element matrix.  We assume that the basis functions and any needed derivatives are already tabulated at the quadrature points.  
After this, the key quantities to consider are the number of degrees of freedom (or local basis functions) of the finite element and the number of quadrature points required.
If we have $N_{dof}$ basis functions and $N_q$ quadrature points, we can expect the cost of forming the element matrix to be proportional to the quantity
\[
C \equiv N_{dof}^2 N_q.
\]
This value of $C$ will be multiplied by the cost of the innermost loop nest, which depends on the particulars of the variational form under consideration.

The particular values for $N_{dof}$ and $N_q$ will depend on the polynomial degree and $|\Delta(K)|$, the number of subcells into which each cell $K$ is split.
Macroelements tend to have a lower $N_{dof}$ than classical ones, but the size of $N_q$ is more complicated.
Although one may use a simpler quadrature rule with lower degree polynomials, one must tile that quadrature rule over each subcell.
We may refine our value of $C$ by noting that
$N_q = N_q^{ref} |\Delta(K)|$, the size of a single-cell quadrature rule times the number of subcells.  Hence, the total work is proportional to
\begin{equation}
  \label{eq:C}
C = N_{dof}^2 N_q^{ref} |\Delta(K)|.
\end{equation}
To further fix ideas, with polynomials of degree $k$, we employ a quadrature rule that is exact on polynomials of degree $2k$.
Specifically,  we use the Xiao-Gimbutas quadrature rules~\cite{xiao2010numerical} now available in FIAT~\cite{brubeck2024fiat}.
Table~\ref{table:C} gives the value of $C$ for several $C^1$ element choices.
For example, PS12 uses only quadratic polynomials and a low degree quadrature rule, having to tile that rule over 12 subcells makes the element matrix almost as expensive to compute as for the much-larger quintic Argyris element.
On the other hand, both PS6 and HCT3 turn out to be much cheaper.
However, as we show in our timing results later, the cost of solving the linear systems significantly larger than matrix assembly, so additional costs at the level of element assembly need not deter application of macroelements.

\begin{table}
  \begin{tabular}{c|ccccc|c}
Element & Degree & $N_{dof}$ & $|\Delta(K)|$ & $N_q^{ref}$ & $N_q$ & $C$ \\ \hline
PS6 & 2 & 9 & 6 & 6 & 36 & 2916\\
PS12 & 2 & 12 & 12 & 6 & 72 & 10368\\
HCT3 & 3 & 12 & 3 & 12 & 36 & 5184\\
HCT4 & 4 & 19 & 3 & 16 & 48 & 17328\\
A5 & 5 & 21 & 1 & 25 & 25 & 11025\\
  \end{tabular}
  \caption{The cost of forming an elementwise matrix over macroelement spaces is driven by the number local basis functions, the size of a reference quadrature rule, and the number of subcells, per~\eqref{eq:C}.  Costs assume exact quadrature for polynomials of twice the given degree.}
  \label{table:C}
\end{table}

An additional cost not considered here is the cost of applying geometrically-dependent transformations to the basis functions and their derivatives, as described below.
Since this cost is proportional only to $N_{dof} N_q$ rather than $N_{dof}^2 N_q$,
our present calculation gives the general idea of some of the relative costs.
At any rate, since the cost of solving linear systems typically dominates that of assembly, considerations other than optimizing the time spent in local computation may very well drive the choice of a particular element.

\section{Code development}
\label{sec:code}

Here, we describe the development within FIAT itself necessary to enable reference element construction of macroelements,
and the development through other packages to enable full integration within Firedrake.

\subsection{FIAT development}
To begin, we introduced a \lstinline{SimplicialComplex} class in FIAT.
This extension of the existing reference simplices allows multiple cells.
The class defines the local topology/connectivity of the subcells and geometric information such as normal and tangent vectors to the facets in the complex.
The \lstinline{SplitSimplicialComplex} and its subclasses encode particular splittings of an existing simplicial complex.
The constructor for this class takes the complex to be split, a list of the locations of vertices, and topology of the new complex from which
essential parent-to-child relationships are programmatically constructed.

We also provide a \lstinline{MacroQuadrature} class that tiles a given simplicial quadrature rule over each subcell or facet to give composite quadrature rules over a splitting.
This works with any of FIAT's wide suite of general quadrature rules,
such as the recently-added Xiao--Gimbutas rules~\cite{xiao2010numerical}.
When finite element methods employ elements over different splittings, one requires quadrature rules accurate on both complexes.  To support this,
simplicial complexes also provide a kind of comparison operator indicating when one complex is a refinement of another.

We maintain FIAT's approach of building basis functions from an \emph{expansion set} -- some set of orthonormal polynomials -- with macroelements.
We can tile $L^2$ orthogonal polynomials~\cite{karniadakis2005spectral}, computed using the recurrences in~\cite{kirby2010singularity}, across the subcomplex without continuity.
We can also tile the modified $C^0$ expansion set described in~\cite{karniadakis2005spectral} continuously over splittings.
Bases for spaces with higher continuity are found as the null space of a collection of functionals defining jumps across internal boundaries.

FIAT constructs the nodal basis for a finite element -- including macroelements -- by pairing a basis for the approximating space with a list of degrees of freedom containing a basis for the dual space, as described in prior work~\cite{Kirby:2004}.
This key aspect of FIAT has required no further internal modification, but
care must be taken that functionals defining integral moments over a split complex use an appropriate composite quadrature rule.

Lagrange finite elements over macro cells may be constructed in one of two ways.  First, the user may instantiate the \lstinline{Lagrange} finite element over a \lstinline{SplitSimplicialComplex}, in which case a suitable expansion set is chosen.  Alternatively, one may provide an unsplit simplex and the optional \lstinline{variant} keyword.  When this keyword provides a splitting such as \lstinline{'iso'} or \lstinline{'alfeld'}, the reference simplex is appropriately split.

\subsection{FInAT development}
FInAT~\cite{homolya2017exposing} provides the main basis function interface to the rest of Firedrake.
It provides abstract syntax for evaluating and manipulating basis functions, enabling optimizations such as sum-factorization, along with many other services.
Critically, FInAT provides a mechanism to implement the non-standard transformations for finite elements with derivative degrees of freedom~\cite{finat-zany}.
Some of the basis transformations are described below.

\subsection{TSFC development}
Firedrake uses \lstinline{tsfc}~\cite{homolya2018tsfc} to generate code for evaluating variational forms.
Implementing macroelements has required only minor changes within \lstinline{tsfc}.
In particular, the algorithm for quadrature selection had to be generalized to select a composite rule suitable for all functions (including those defined over splittings) appearing in a given integral in the form.

Macroelements present an opportunity for future optimizations within
the form compiler.  Currently, \lstinline{tsfc} uses a ``flattened''
quadrature rule and table of basis values that has no reference to
subelement structure.
Alternatively, \lstinline{tsfc} could handle macroelements by spawning an inner loop over the subcells, assembling the element-level bases from subcell-level shape functions.
With high-continuity elements, one typically has that most local basis functions are supported on each subcell in a split so that we might not expect a significant speedup from this approach.
On the other hand, $C^0$ elements defined over higher-order iso-type
splits could benefit since the basis functions would be locally
supported on only a few of the subcells.
Here, we have focused on feature inclusion rather than the
potentially invasive modifications required to optimize this
particular case.

\subsection{Firedrake development}
Somewhat surprisingly, our design of macroelements in FIAT and its interfaces with FInAT and \lstinline{tsfc} are rather self-contained and spawned 
no internal Firedrake development.  One simply obtains a macroelement space by providing the element name to the \lstinline{FunctionSpace} constructor.
To obtain Lagrange-type macroelements, one provides the \lstinline{variant} keyword.  For example:
\begin{lstlisting}
  V0 = FunctionSpace(mesh, 'HCT', 3)
  V1 = FunctionSpace(mesh, 'Lagrange', 2, variant='alfeld')
\end{lstlisting}
In fact, the systematic modifications in the supporting packages meant that the Firedrake pull request only added tests of accuracy with macroelements.

\section{Transformation theory}
\label{sec:xform}
\begin{figure}[htbp]
  \centering
  \begin{tikzpicture}
    \draw[fill=lightgray] (0,0) coordinate (vhat1)
    -- (2,0) coordinate(vhat2)
    -- (0,2) coordinate(vhat3)--cycle;
    \foreach \pt\labpos\lab in {vhat1/below/\hat{\mathbf{v}}_1, vhat2/below right/\hat{\mathbf{v}}_2, vhat3/above right/\hat{\mathbf{v}}_3}{
      \filldraw (\pt) circle(.6mm) node[\labpos=1.5mm, fill=white]{$\lab$};
    }
    %\draw[->] (vhat2) -- (2.5, 0);
    %\draw[->] (vhat3) -- (0, 2.5);
    \draw[fill=lightgray] (4.0, 1.5) coordinate (v1)
    -- (5.5, 2.0) coordinate (v2)
    -- (4.8, 2.7) coordinate (v3) -- cycle;
    \foreach \pt\labpos\lab in {v1/below/\mathbf{v}_1, v2/below right/\mathbf{v}_2, v3/above right/\mathbf{v}_3}{
      \filldraw (\pt) circle(.6mm) node[\labpos=1.5mm, fill=white]{$\lab$};
    }
    \draw[<-] (1.1, 1.1) to[bend left] (4.3, 2.2);
    \node at (2.7, 2.65) {$F:K\rightarrow\hat{K}$};
    \node at (0.6,0.6) {$\hat{K}$};
    \node at (4.75, 2.1) {$K$};
  \end{tikzpicture}
  \caption{Affine mapping to a reference cell \(\hat{K}\) from a
    typical cell \( K \).  Note that here $F$ maps from the physical
    cell $K$ to the reference cell $\hat{K}$ rather than the other way
    around.}
  \Description{Diagram showing affine map of a physical triangle to the reference triangle.}
    \label{fig:affmap}
\end{figure}
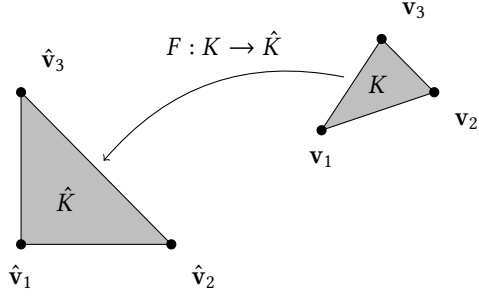

When basis functions are constructed on some reference domain $\hat{K}$, they must be mapped to each cell $K$ in the given mesh.
Classically, a mapping $F: K \rightarrow \hat{K}$ (affine for straight-sided simplices) is constructed per ~\Cref{fig:affmap}, and one obtains a pullback in the standard way.
For each
$\hat{f} : \hat{K} \rightarrow \mathbb{R}$, define $F^*(\hat{f}) \equiv f : K \rightarrow \mathbb{R}$ by
\begin{equation}
  f = F^*(\hat{f}) = \hat{f} \circ F,
\end{equation}
and we also have the push-forward of a functional $n$ acting on functions over $K$ to those over $\hat{K}$ by
\begin{equation}
  F_*(n) = n \circ F^*.
\end{equation}
Note that, as in~\cite{BreSco} we are defining $F$ to map from the physical cell to the reference cell.
Many references take $F$ going in the reverse direction (which may be more natural as this is more directly computable in the non-affine case), but having it go from physical to reference cell obviates the need for inverses in the notation for pullback and push-forward.

\begin{figure}[htbp]
  \centering
  \begin{tikzpicture}
    \draw[fill=magenta] (0,0) -- (2,0) -- (0,2) -- cycle;
    \foreach \x\y in {0/0, 2/0, 0/2}{
      \filldraw (\x,\y) circle(0.04);
      \draw[->] (\x,\y) -- (\x+6/7*0.2,\y-4/7*0.2);
      \draw[->] (\x,\y) -- (\x-5/14*0.2, \y+15/14*0.2);
    }
    \draw[dashed] (0,0) -- (2/3,2/3);
    \draw[dashed] (2,0) -- (2/3,2/3);
    \draw[dashed] (0,2) -- (2/3,2/3);
    \foreach \x\y\nx\ny in {1/1/.3536/.3536,0/1/-.9113/1.070,1/0/.6010/-1.1971}{
      \draw[->] (\x, \y) -- (\x+0.2*\nx, \y+0.2*\ny);
      }
    
    \draw[fill=magenta] (4.0, 1.5) -- (5.5, 2.0) -- (4.8, 2.7) -- cycle;
    \foreach \x\y in {4/1.5, 5.5/2, 4.8/2.7}{
      \filldraw (\x,\y) circle(.04);
      \draw[->] (\x,\y)--(\x,\y+0.2);
      \draw[->] (\x,\y)--(\x+0.2,\y);
    }
    \draw[dashed] (4,1.5) -- (4.77,2.07);
    \draw[dashed] (5.5,2) -- (4.77,2.07);
    \draw[dashed] (4.8,2.7) -- (4.77,2.07);
    \foreach \x\y\nx\ny in {5.15/2.35/.707/.707,4.4/2.1/-0.8321/0.5547,4.75/1.75/.3163/-.9487}{
      \draw[->] (\x,\y) -- (\x+0.2*\nx,\y+0.2*\ny);
      }

    \draw[<-] (1.2, 1.1) to[bend left] (4.2, 2.2);
    \node at (2.7, 2.65) {$F_*$};
    \end{tikzpicture}
  \caption{Pushing forward the HCT derivative nodes in physical
    space does \emph{not} produce the reference derivative nodes.}
  \Description{Diagram showing the failure of affine push-forward to preserve nodes in the HCT element.}
\label{fig:hctpushforward}
\end{figure}
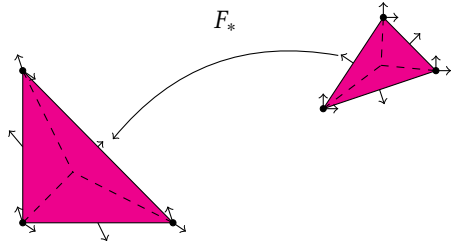

For Lagrange elements, the pullback maps the reference element basis functions exactly to the physical element basis functions, and the push-forward maps point evaluation on $K$ to point evaluation on the point's preimage points under $F$.
However, $C^1$ and many other finite elements utilize derivatives and other degrees of freedom that are not preserved under push-forward, as shown in~\Cref{fig:hctpushforward}, and this complicates the mappings.
Similar issues are observed for $\hdiv$ finite elements and the contravariant Piola map when one has degrees of freedom involving tangents or vertex values~\citep{aznaran2022transformations}.
Here, we apply the theory developed in~\citep{kirby-zany} to some of our newly-enabled macroelements.
This theory lets us  identify a matrix $M$ such that the vector of basis functions over $K$ can be constructed by $M$ times the vector of pullbacks of reference element basis functions.
The transpose of this matrix, relating the push-forwards of physical nodes to the reference element nodes, is often easier to construct. 
This matrix is typically quite sparse, making this approach much cheaper than directly building the basis functions on each cell.

We let $\mathcal{N}$ be a vector containing the finite element nodes over the physical element and $\hat{\mathcal{N}}$ be the reference element nodes.
Then, if $\mathcal{F}_*(\mathcal{N})$ is defined componentwise as the push forward of each entry of $\mathcal{N}$, we seek a matrix $V$ such that
\begin{equation}
  \hat{\mathcal{N}} = V \mathcal{F}_*(\mathcal{N}).
\end{equation}
As shown in~\cite{kirby-zany}, we have that $M = V^\top$.

When the function space is preserved under pullback and we have affine-interpolation equivalent elements (like the Hermite simplex), $V$ has a relatively simple block-diagonal structure.
For other elements such as Morley and Argyris, the situation is more complex.
We introduce a ``completion'' of the nodes, contained in vectors of functionals $\mathcal{N}^c$ and $\hat{\mathcal{N}}^c$.
These new vectors include all entries present in $\mathcal{N}$ and $\hat{\mathcal{N}}$ and have the property that $\operatorname{span}(F_*(N^c)) = \operatorname{span}(\hat{\mathcal{N}}^c)$ in $C^k(\hat{K})^\prime$.
This choice is often apparent from context, as we will see below for the HCT element.  Then, the theory in~\cite{kirby-zany} shows that we can obtain $V = M^\top$ in factored form by
\begin{equation}
  \label{eq:factoredV}
  V = E V^c D.
\end{equation}
Here, the matrix $D$ maps the physical nodes $\mathcal{N}$ (or at least their restriction to elements of $P$) to their completion $\mathcal{N}^c$.
Then, $V^c$ is the matrix such that
\begin{equation}
  \hat{\mathcal{N}}^c = V^c \mathcal{F}_*(\mathcal{N}^c).
\end{equation}
That is, it relates the push-forwards of the completion of the physical nodes to their reference counterparts.
Finally, the matrix $E$ simply extracts the members $\hat{\mathcal{N}}^c$
that also belong in $\hat{\mathcal{N}}$ -- it has only one nonzero entry per row.

\subsection{HCT}
The $C^1$-conforming HCT triangle has 12 degrees of freedom and piecewise cubic polynomials.
Each $K$ is divided into three triangles $K_1$, $K_2$, $K_3$ by the Alfeld split.
The function space $P^{HCT}(K)$ then consists of cubic polynomials over each $K_i$ that are $C^1$ across internal edges.
This space is parameterized by function values and gradients at the vertices together with normal derivatives on each edge.

%% To fix ideas, we denote point evaluation at some point $\bfx$ by
%% \begin{equation}
%%   \delta_\bfx(p) = p(\bfx).
%% \end{equation}
%% We let $\delta_\bfx^\bfx$ denote the derivative in some
%% direction $\bfs$ at a point $\bfx$:
%% \begin{equation}
%%     \delta^\bfs_\bfx(p) = \bfs\cdot \nabla p(\bfx)
%% \end{equation}

%% \begin{equation}
%%   \label{eq:edgederivmoment}
%%   \mu^\bfs_\bfe(f) = \int_\bfe \bfs \cdot \nabla f \, ds, \\
%% \end{equation}

%% We use block notation, with
%% \begin{equation}
%%   \nabla_\bfx = \begin{bmatrix} \delta_\bfx^\bfx &
%%     \delta_{\bfx}^\bfy \end{bmatrix}^\top
%% \end{equation}
%% for the gradient in Cartesian coordinates at a point $\bfx$.

We can write the nodes for HCT in a vector as
\begin{equation}
  \label{eq:hctnodes}
  \mathcal{N} =
  \begin{bmatrix}
    \delta_{\bfv_1} & \nabla^\top_{\bfv_1} &
    \delta_{\bfv_2} & \nabla^\top_{\bfv_2} &
    \delta_{\bfv_3} & \nabla^\top_{\bfv_3} &
    \mu^{\bfn_1}_{\bfe_1} & \mu^{\bfn_2}_{\bfe_2} &
    \mu^{\bfn_3}_{\bfe_3}
  \end{bmatrix}^\top.
\end{equation}

For the reference element, we define the edge nodes to use integral averages rather than moments:
\begin{equation} \label{eq:edgederivmomentref}
  \hat{\mu}^{\hat{\bfs}}_{\hat{\bfe}}(f) = \tfrac{1}{|\hat{\bfe}|}\int_{\hat{\bfe}} \hat{\bfs} \cdot \hat{\nabla} f \, d\hat{s}.
\end{equation}
Then, we enumerate the reference element nodes as
\begin{equation}
  \label{eq:hctrefnodes}
  \widehat{\mathcal{N}} =
  \begin{bmatrix}
    \delta_{\hat{\bfv}_1} & \hat{\nabla}^\top_{\hat{\bfv}_1} &
    \delta_{\hat{\bfv}_2} & \hat{\nabla}^\top_{\hat{\bfv}_2} &
    \delta_{\hat{\bfv}_3} & \hat{\nabla}^\top_{\hat{\bfv}_3} &
    \hat{\mu}^{\hat{\bfn}_1}_{\hat{\bfe}_1} &
    \hat{\mu}^{\hat{\bfn}_2}_{\hat{\bfe}_2} &
    \hat{\mu}^{\hat{\bfn}_3}_{\hat{\bfe}_3}
  \end{bmatrix}^\top.
\end{equation}
This redefinition eliminates the need for logic indicating which edges of each triangle correspond to which reference element edges.

Now, for the HCT element, the push-forward of the physical nodes does not align with the reference nodes, and so we need the theory outlined above.
To form the completion, we first introduce $\mu^{\bft}_\bfe(f)$ per~\eqref{eq:edgederivmoment} to be the integral moment of the tangential derivative along edge $\bfe$ of the triangle.
To find the nodal completion for the reference dual, we take $\hat{\mu}^{\hat{\bft}}_{\hat{\bfe}}$ to be the integral average of the tangential derivative along a reference edge.
We define
\begin{equation}
  \mathcal{M}_{i} =
  \begin{bmatrix}
    \mu^{\bfn_i}_{\bfe_i} &
    \mu^{\bft_i}_{\bfe_i} \end{bmatrix}^\top
\end{equation}
to be the vector of the moments of the normal and tangential
derivatives on a particular edge.
We also let $\widehat{\mathcal{M}_{i}}$ contain the corresponding reference element nodes.

Then, the compatible nodal completion for the HCT element has physical nodes
\begin{equation}
  \label{eq:hctnodescompleted}
  \mathcal{N}^c =
  \begin{bmatrix}
    \delta_{\bfv_1} & \nabla^\top_{\bfv_1} &
    \delta_{\bfv_2} & \nabla^\top_{\bfv_2} &
    \delta_{\bfv_3} & \nabla^\top_{\bfv_3} &
    \mathcal{M}_1^\top &
    \mathcal{M}_2^\top &
    \mathcal{M}_3^\top
  \end{bmatrix}^\top.
\end{equation}
and reference nodes
\begin{equation}
  \label{eq:hctrefnodescompleted}
  \widehat{\mathcal{N}}^c =
  \begin{bmatrix}
    \delta_{\hat{\bfv}_1} & \hat{\nabla}^\top_{\hat{\bfv}_1} &
    \delta_{\hat{\bfv}_2} & \hat{\nabla}^\top_{\hat{\bfv}_2} &
    \delta_{\hat{\bfv}_3} & \hat{\nabla}^\top_{\hat{\bfv}_3} &
    \widehat{\mathcal{M}_1}^\top &
    \widehat{\mathcal{M}_2}^\top &
    \widehat{\mathcal{M}_3}^\top
  \end{bmatrix}^\top.
\end{equation}

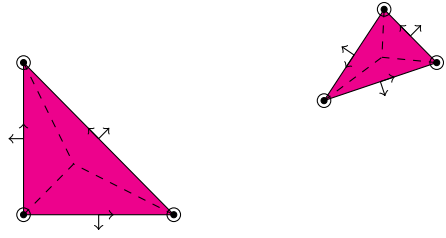
\begin{figure}[htbp]
  \centering
  \begin{tikzpicture}
    \draw[fill=magenta] (0,0) -- (2,0) -- (0,2) -- cycle;
    \foreach \x\y in {0/0, 2/0, 0/2}{
      \filldraw (\x,\y) circle (0.04);
      \draw (\x, \y) circle (0.09);
    }
    \draw[dashed] (0,0) -- (2/3,2/3);
    \draw[dashed] (2,0) -- (2/3,2/3);
    \draw[dashed] (0,2) -- (2/3,2/3);    

    \foreach \x\y\nx\ny\tx\ty in {1/1/.7071/.7071/-.7071/.7071,
      0/1/-1/0/0/1,
      1/0/0/-1/1/0}
    {
      \draw[->] (\x, \y) -- (\x+0.2*\nx, \y+0.2*\ny);
      \draw[->] (\x, \y) -- (\x+0.2*\tx, \y+0.2*\ty);
    }
    
    \draw[fill=magenta] (4.0, 1.5) -- (5.5, 2.0) -- (4.8, 2.7) -- cycle;
    \foreach \x\y in {4/1.5, 5.5/2, 4.8/2.7}{
      \filldraw (\x,\y) circle(.04);
      \draw (\x, \y) circle (0.09);      
    }
    \draw[dashed] (4,1.5) -- (4.77,2.07);
    \draw[dashed] (5.5,2) -- (4.77,2.07);
    \draw[dashed] (4.8,2.7) -- (4.77,2.07);    
    \foreach \x\y\nx\ny\tx\ty in {5.15/2.35/.707/.707/-.707/.707,4.4/2.1/-0.8321/0.5547/-.5547/-.8321,4.75/1.75/.3163/-.9487/.9487/.3163}{
      \draw[->] (\x,\y) -- (\x+0.2*\nx,\y+0.2*\ny);
      \draw[->] (\x,\y) -- (\x+0.2*\tx,\y+0.2*\ty);
      }
    \end{tikzpicture}
  \caption{Nodal sets \( \hat{N}^c \) and \( N^c \) giving the
    compatible nodal completion of $N$ and $\hat{N}$ for an HCT 
    element and reference element are formed by including tangential
    derivatives along with normal derivatives at each edge midpoint.}
  \Description{Diagram showing the nodal completion for HCT elements.}
\label{fig:hctbridge}
\end{figure}

The space is 12-dimensional, so the lengths of $\mathcal{N}$ and
$\hat{\mathcal{N}}$ are both 12.
The completions contain three extra entries, so $\mathcal{N}^c$ and $\hat{\mathcal{N}}^c$ have length 15.
So, the factors from~\eqref{eq:factoredV} have
\begin{equation}
  \begin{split}
    E & \in \mathbb{R}^{12 \times 15}, \\
    V^c & \in \mathbb{R}^{15 \times 15}, \\
    D & \in \mathbb{R}^{15 \times 12}.
  \end{split}
\end{equation}

Now, we specify these matrices for HCT from right to left (as they are applied to a vector).
The form of $D$ follows quite naturally.
Clearly, the rows corresponding to members of $\mathcal{N}^c$ also appearing in $\mathcal{N}$ will just have a single nonzero in the appropriate column.  
The nodes in $\mathcal{N}^c$ not in $\mathcal{N}$ are just integrals of quantities over edges, and we can use the Fundamental Theorem of Calculus to perform this task.
Let $\bfe$ be an edge running from vertex $\bfv_a$ to $\bfv_b$ with unit tangent $\bft$.  We have
\begin{equation}
\label{eq:mut}
  \mu^{\bft}_{\bfe}(f) = \int_{\bfe} \bft\cdot \nabla f \ds
  = f(\bfv_b)-f(\bfv_a) = \delta_{\bfv_b}(f) - \delta_{\bfv_a}(f)
\end{equation}

Putting this together with the node orderings in~\eqref{eq:hctnodescompleted} and~\eqref{eq:hctnodes}, we write $D$ as a block matrix segregating the vertex and edge nodes with
%% \begin{equation}
%%   D = \begin{bmatrix}
%%     1 & 0 & 0 & 0 & 0 & 0 & 0 & 0 & 0 & 0 & 0 & 0 \\
%%     0 & 1 & 0 & 0 & 0 & 0 & 0 & 0 & 0 & 0 & 0 & 0 \\
%%     0 & 0 & 1 & 0 & 0 & 0 & 0 & 0 & 0 & 0 & 0 & 0 \\
%%     0 & 0 & 0 & 1 & 0 & 0 & 0 & 0 & 0 & 0 & 0 & 0 \\
%%     0 & 0 & 0 & 0 & 1 & 0 & 0 & 0 & 0 & 0 & 0 & 0 \\
%%     0 & 0 & 0 & 0 & 0 & 1 & 0 & 0 & 0 & 0 & 0 & 0 \\
%%     0 & 0 & 0 & 0 & 0 & 0 & 1 & 0 & 0 & 0 & 0 & 0 \\
%%     0 & 0 & 0 & 0 & 0 & 0 & 0 & 1 & 0 & 0 & 0 & 0 \\
%%     0 & 0 & 0 & 0 & 0 & 0 & 0 & 0 & 1 & 0 & 0 & 0 \\
%%     0 & 0 & 0 & 0 & 0 & 0 & 0 & 0 & 0 & 1 & 0 & 0 \\
%%     0 & 0 & 0 & -1 & 0 & 0 & 1 & 0 & 0 & 0 & 0 & 0 \\
%%     0 & 0 & 0 & 0 & 0 & 0 & 0 & 0 & 0 & 0 & 1 & 0 \\
%%     -1 & 0 & 0 & 1 & 0 & 0 & 0 & 0 & 0 & 0 & 0 & 0 \\
%%     0 & 0 & 0 & 0 & 0 & 0 & 0 & 0 & 0 & 0 & 0 & 1 \\
%%     -1 & 0 & 0 & 0 & 0 & 0 & 1 & 0 & 0 & 0 & 0 & 0 
%%     \end{bmatrix}
%% \end{equation}
\begin{equation}
  D = \left[ \begin{array}{c|c}
    I & 0 \\ \hline
    D_{12} & D_{22}
    \end{array} \right],
\end{equation}
where, $I$ is the $9 \times 9$ identity matrix, and
\begin{equation}
  D_{12} =
  \begin{bmatrix}
    0 & 0 & 0 & 0 & 0 & 0 & 0 & 0 & 0  \\
    0 & 0 & 0 & -1 & 0 & 0 & 1 & 0 & 0 \\
    0 & 0 & 0 & 0 & 0 & 0 & 0 & 0 & 0 \\
    -1 & 0 & 0 & 1 & 0 & 0 & 0 & 0 & 0 \\
    0 & 0 & 0 & 0 & 0 & 0 & 0 & 0 & 0 \\
    -1 & 0 & 0 & 0 & 0 & 0 & 1 & 0 & 0\\
  \end{bmatrix},  \ \ \
  D_{22} = 
  \begin{bmatrix}
    1 & 0 & 0 \\
    0 & 0 & 0 \\
    0 & 1 & 0 \\
    0 & 0 & 0 \\
    0 & 0 & 1 \\
    0 & 0 & 0
  \end{bmatrix}.
\end{equation}

Then, $V^c$ is obtained by relating the push-forwards of $\mathcal{N}^c$ to $\widehat{\mathcal{N}}^c$.
The vertex degrees of freedom transform just as those for the Hermite triangle~\cite{kirby-zany} and lead to blocks with 1 and the cell Jacobian.
To transform the edge degrees of freedom, we note that
$\mu_\bfe^\bfn$ and $\mu^\bft_e$ together comprise the moments of the gradient along edge $e$ in an orthogonal coordinate system.
Given the orthogonal unit vectors $\bfn$ and $\bft$, we can define an
orthogonal matrix $G$ by:
\begin{equation}
  G = \begin{bmatrix} \bfn & \bft \end{bmatrix}^\top.
\end{equation}
Define $G_i$ to have the normal and tangential
vectors to edge $i$ of triangle $K$ in its columns and $\widehat{G}_i$ those for
triangle $\widehat{K}$. By changing the gradient to normal/tangential coordinates,
\begin{equation}
  \label{eq:gradtont}
  \nabla_x = G^\top \nabla^{\bfn\bft}_\bfx.
\end{equation}
Physical and reference gradients in the normal/tangential coordinates follow via the chain rule
\begin{equation}
  \label{eq:transformnortangrad}
  \nabla_{\bfx}^{\bfn\bft} = G J^\top \widehat{G}^\top \hat{\nabla}^{\hat{\bfn}\hat{\bft}}_{\hat{\bfx}}.
\end{equation}
Now, for any vector $\bfs$, edge $\bfe$, and smooth function $f = f \circ F$, we have
\begin{equation}
\int_\bfe \bfs\cdot \nabla f \ds =
\int_\bfe \bfs\cdot \hat\nabla f \circ F \ds 
= \int_{\hat{\bfe}} \bfs\cdot \hat\nabla f J_{\bfe,\hat{\bfe}} d\hat{s},
\end{equation}
where the Jacobian $J_{\bfe,\hat{\bfe}}$ is just the ratio of the
length of $\bfe$ to that of the corresponding reference element edge
$\hat{\bfe}$. Applying this to the normal and tangential moments and
using~\eqref{eq:transformnortangrad}, we have that:
\begin{equation}
\mathcal{M}_{i}
= |\bfe_i| G_i J^\top \hat{G}_{i}^\top \widehat{\mathcal{M}}_{i},
\end{equation}
where the factor of $|\hat{\bfe}_i|$ in the denominator of the
Jacobian is merged with the reference element moments to produce
$\widehat{\mathcal{M}_{i}}$.  Hence, the slight modification of
reference element nodes avoids extra data structures or logic in
identifying reference element edge numbers.

We identify the matrices
\begin{equation}
  B_{i} = \hat{G}_i^{\top} J^{-\top} G / |\bfe_i|.
\end{equation}
Then, we can also write $V^c$ in a block diagonal form
\begin{equation}
  V^c = \left[ \begin{array}{c|c}
    V^c_v & 0 \\ \hline
    0 & V^c_e 
    \end{array} \right],
\end{equation}
where $V^c_v$ is itself a block $3 \times 3$ matrix with three copies of
the $3 \times 3$ matrix
\begin{equation}
\begin{bmatrix} 1 & 0 \\
  0 & J^{-T}
  \end{bmatrix}
\end{equation}
along the diagonal.  $V^c_e$ is also a block $3 \times 3$ matrix
\begin{equation}
V^c_e = \begin{bmatrix} B_1 & 0 & 0 \\ 0 & B_2 & 0 \\ 0 & 0 & B_3 \end{bmatrix}.
\end{equation}

The Boolean extraction matrix $E \in \mathbb{R}^{12 \times 15}$ simply takes entries of $\widehat{N}$ from $\widehat{\mathcal{N}}^c$, so that
\[
E_{ij} = \begin{cases} 1, & i = j \text{ and } 1 \leq i \leq 9 \text{ or } i,j \in \left\{ (10, 10), (11, 12), (12, 14) \right\}, \\
  0, & \text{ otherwise}.
  \end{cases}
\]

The reduced HCT has fewer degrees of freedom but requires a
more complicated transformation.  The affine map takes functions with a
degree-constrained normal derivative in reference space to
functions with a degree constraint in some other direction.
The Bell element (quintics with cubic normal
derivatives on each edge) provides a classical example of this feature in a non-macroelement.
A generalized theorem from~\cite{kirby-zany} handles such cases. 
Essentially, one interprets the constraints as additional degrees of freedom in an
``extended'' finite element over the unconstrained space.
This element can be transformed by the above theory, and the correct basis functions for the reduced space selected afterward.
The adaptation of the HCT transformation to the reduced case is
quite analogous to constructing the Bell transformation from the
Argyris in~\cite{kirby-zany}.

\subsection{High-order HCT}
We briefly describe the generalization to higher-order HCT, using a similar
space as the one from \cite{groselj22}.  In an effort to reduce the number
of interior degrees of freedom while maintaining the same order of
accuracy as the minimally $C^1$ polynomials of degree $k \ge 3$,
they enforce $C^{k-1}$ supersmoothness at the interior vertex of the
split.
To enable use of a reference element, we take this interior vertex as the barycenter,
although the incenter could be used.  We obtain the
$C^1$ subspace by constraining the $C^0$ hierarchical basis on an Alfeld split
using the barycenter.

The high-order HCT nodes are an augmentation of \eqref{eq:hctnodes}.
The vertex-based degrees of freedom remain unchanged, i.e.,
we keep point evaluation and the gradient at each exterior vertex.

The edge-based degrees of freedom come in two kinds.
The first kind are normal derivative moments against Jacobi polynomials $P_i^{(1,1)}$,
\begin{equation}
   \hat{\mu}^{\bfn}_{\bfe,i}(f) = 
   \tfrac{1}{|\bfe|} \int_\bfe P_i^{(1,1)}(\hat{s}) \bfn\cdot\nabla f \ds,
   \quad i=0,\ldots,k-3.
\end{equation}
Here, we introduce the normalized variable $\hat{s}\in [-1, 1]$, 
defined in terms of the length $s$ along $\bfe$ as $\hat{s} = 2(s/|\bfe|) - 1$. 
We have slightly altered the notation, using the second subscript to index into a set of polynomials
rather than as a polynomial itself.
The second kind are trace moments against the derivative of Jacobi polynomials 
\begin{equation}
   \mu_{\bfe,i}(f) = 
   \int_\bfe \frac{d}{ds}P_i^{(1,1)}(\hat{s}) f \ds,
   \quad i=1,\ldots,k-3.
\end{equation}
For $k \ge 4$, the interior degrees of freedom
are integral moments against a basis $\{q_i\}$ for $\P_{k-4}(K)$,
\begin{equation}
   \mu_{K,i}(f) = \int_K q_i f \dx. 
\end{equation}

For the edge degrees of freedom, the Jacobi weights $(1,1)$ have been 
chosen to produce hierarchical edge-based basis functions. Using the orthogonality property
\begin{equation}
   \int_{-1}^1 (\hat{s}+1)^\alpha (\hat{s}-1)^\beta P_i^{(\alpha,\beta)}(\hat{s}) P_j^{(\alpha,\beta)}(\hat{s}) \, \mathrm{d}\hat{s}
   = c^{(\alpha,\beta)}_i \delta_{ij},
\end{equation}
the identity $\tfrac{d}{d\hat{s}} P_i^{(\alpha,\beta)}(\hat{s}) = d^{(\alpha,\beta)}_i P_{i-1}^{(\alpha+1,\beta+1)}(\hat{s})$, 
and the enforcement of vanishing value and tangential derivative at the vertices,
we deduce that the edge-based basis functions associated with $\mu_{\bfe,i}$ have a trace
proportional to $(\hat{s}+1)^2(\hat{s}-1)^2 P_{i-1}^{(2,2)}(\hat{s})$.
Similarly, basis functions associated with $\mu^{\bfn}_{\bfe,i}$
have normal derivative trace proportional to $(\hat{s}+1)(\hat{s}-1) P_i^{(1,1)}(\hat{s})$.

On each edge $\bfe$, for a fixed value of $i > 0$, to find the nodal completion
of the normal derivative moments $\mu^\bfn_{\bfe,i}$, 
again we take the tangential derivative moments $\mu^\bft_{\bfe,i}$. 
These can conveniently computed from the trace
moments $\mu_{\bfe,i}$ and the values of $f$ at the endpoints $\bfv_a, \bfv_b$ of $\bfe$, by
integrating by parts and exploiting the fact that the trace moments are defined against $\tfrac{\mathrm{d}}{\mathrm{d}s} P^{(1,1)}_i$,
\begin{equation}
   \mu^{\bft}_{\bfe,i}(f) =
   \int_\bfe  P_i^{(1,1)}(\hat{s}) \bft\cdot\nabla f \ds = 
   - \mu_{\bfe,i}(f) + P_i^{(1,1)}(1) \delta_{\bfv_b}(f) - P_i^{(1,1)}(-1)\delta_{\bfv_a}(f).
\end{equation}
Note that \eqref{eq:mut} still holds for $i=0$. 
The transformation of $\mu^\bfn_{\bfe,i}$ then proceeds in a very similar fashion,
except that there is additional coupling with $\mu_{\bfe,i}$. 

We may also extend this construction to Argyris elements of arbitrarily high order.
The only difference is that the Jacobi weights need to be chosen as $(2,2)$,
due to the presence of second derivative degrees of freedom at vertices.

\section{Numerical examples}
\label{sec:num}

Our computations were performed with Firedrake installed
on a Linux workstation with an AMD Ryzen Threadripper PRO 7975WX CPU
with 32 cores, 64 threads, a clock-rate of 2.194 GHz, and 768 GiB of RAM.

\subsection{Stokes flow}
Here, we consider two formulations of the Stokes equations of incompressible flow for which our our macroelement technology provides effective discretization.
First, we consider the standard pressure-velocity formulation, written in weak form as
\begin{equation}
  \label{eq:velprestokes}
  \begin{split}
     \left( 2\nu\upvarepsilon(\mathbf{u}), \upvarepsilon(\mathbf{v}) \right) - \left(p, \div\mathbf{v}\right) & = \left(\mathbf{f}, \mathbf{v}\right), \\
    -\left( \div\mathbf{u}, q \right) & = 0.
  \end{split}
\end{equation}
Here, $\mathbf{u}$ and $p$ represent the fluid velocity and pressure in some domain $\Omega$ in two or three spatial dimensions.
The parameter $\nu$ is known as the kinematic viscosity,
and the system may be driven by body forces $\mathbf{f}$ (e.g. gravity), 
and boundary conditions on the velocity and/or stress $\sigma \coloneqq 2\nu \upvarepsilon(\mathbf{u}) - p I$
must be included to close the system. 

The Stokes equations may also be formulated in terms of the symmetric stress
$\sigma$ itself and velocity $\mathbf{u}$, which may be of interest in extensions to  non-Newtonian flows. 
Here we use $\tr\sigma = 2\nu \div\mathbf{u} - d p$
to eliminate the pressure from the system.
We integrate the implicit constitutive relation 
$\upvarepsilon(\mathbf{u}) = \frac{1}{2\nu}(\sigma - \frac{1}{d} \tr{\sigma} I)$
against a symmetric tensor $\tau$, and integrate by parts
$(\upvarepsilon(\mathbf{u}), \tau) = (\mathbf{u}, \tau\mathbf{n})_{\partial\Omega} -(\mathbf{u}, \div\tau)$.
Together with the balance of momentum equation, we arrive at the weak formulation
\begin{equation}
  \label{eq:stressvelstokes}
  \begin{split}
   \frac{1}{2\nu} \left(\sigma, \tau\right) - \frac{1}{2d\nu} \left(\tr\sigma, \tr\tau\right) + \left(\mathbf{u}, \div\tau\right) & = 
     (\mathbf{u}_0, \tau\mathbf{n})_{\partial\Omega},\\
   \left(\div\sigma, \mathbf{v}\right) & = -\left(\mathbf{f}, \mathbf{v}\right).
  \end{split}
\end{equation}
Normal stress boundary conditions may be strongly imposed on part of the boundary,
and on the complement of the boundary, $\mathbf{u}_0$ is the displacement to be weakly imposed.

Discretizing symmetric stress tensors requires special treatment.
One may use the Arnold--Winther elements as in~\cite{carstensen2012numerical}, 
or weakly enforce the symmetry of the stress tensor~\cite{arnold1984peers,gopalakrishnan2020mass}.
The Johnson--Mercier macroelement provides a conforming and symmetric stress
approximation with only piecewise linear polynomials.

Here, we attempt to assess the accuracy and cost of various methods for both
the velocity-pressure and stress-velocity forms of the Stokes equations in two
and three dimensions.  We use the method of manufactured solutions, selecting
the right-hand side and boundary data such that the true solution to the
equations is some known smooth function.  We force both forms of the equation
(velocity-pressure and stress-velocity) to yield the same pressure, velocity,
and stress.

We partition the unit square into an $N \times N$ mesh of squares, each split into two right triangles.
The unit cube is partitioned into an $N \times N \times N$ mesh of cubes, each split into six tetrahedra.
We approximate~\eqref{eq:velprestokes} with several finite element families.
We use the standard $\P_2-\P_1$ Taylor--Hood and the $\P_2^\mathrm{iso}-\P_1$ elements to provide a point of comparison to classical elements.
We also use Scott--Vogelius macroelements with continuous $\P_d$ velocity and discontinuous $\P_{d-1}$ pressures on the Alfeld split,
the Alfeld--Sorokina velocity element paired with $C^0 \P_1$ velocity on the Alfeld split and the Guzm\'an--Neilan element paired with piecewise constants.
Finally, we approximate~\eqref{eq:stressvelstokes} using the Johnson--Mercier macroelement paired with discontinuous linear velocities.

We measure the enforcement of incompressibility constraint $\div\mathbf{u}=0$
with the functional
\begin{equation} \label{eq:divnorm}
   \left(\sum_{K} \| \div\mathbf{u} \|^2_K
   +  \sum_{e} \tfrac{1}{h_e} \| \llbracket \mathbf{u} \cdot \mathbf{n}  \rrbracket \|_e^2\right)^{1/2},
\end{equation}
where $\| \cdot \|_K$ denotes the $L^2$ norm over a cell and $\| \cdot \|_e$ the $L^2$ norm over an edge.  The quantity $\llbracket \cdot \rrbracket$ indicates the jump, or difference in a quantity across an edge of the triangulation.
The velocity-pressure formulation gives an $H^1$-conforming discretization, so this quantity simplifies to the $L^2$ norm of the divergence.

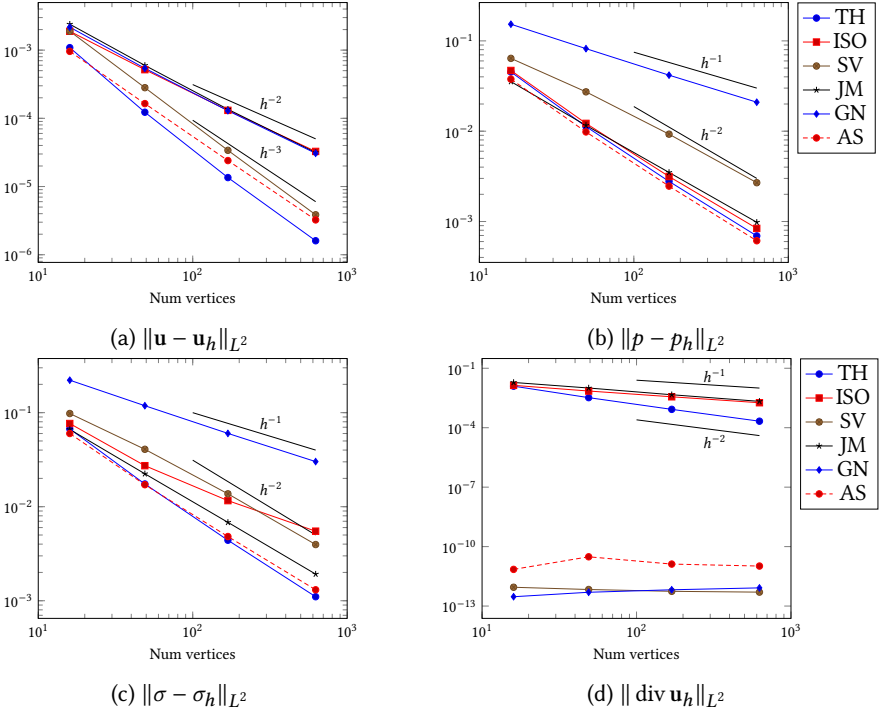
\begin{figure}[htbp]
  \begin{subfigure}[t]{0.45\textwidth} \centering
    \begin{tikzpicture}[scale=0.6]
      \begin{loglogaxis}[xlabel=Num DOFs, legend style={font=\LARGE},
          tick label style={font=\Large},
          label style={font=\Large},
          legend pos=outer north east
         ]
        \addplot table[x=dofs, y=velocityL2, col sep=comma]{code/stokes/data/CG2xCG1.2d.csv}; %\addlegendentry{TH};
        \addplot table[x=dofs, y=velocityL2, col sep=comma]{code/stokes/data/CG1isoxCG1.2d.csv}; %\addlegendentry{ISO};
        \addplot table[x=dofs, y=velocityL2, col sep=comma]{code/stokes/data/CG2xDG1alfeld.2d.csv}; %\addlegendentry{SV};
        \addplot table[x=dofs, y=velocityL2, col sep=comma]{code/stokes/data/JM1xDG1.2d.csv}; %\addlegendentry{JM};
        \addplot table[x=dofs, y=velocityL2, col sep=comma]{code/stokes/data/GN1xDG0.2d.csv}; %\addlegendentry{GN};
        \addplot table[x=dofs, y=velocityL2, col sep=comma]{code/stokes/data/AS2xCG1.2d.csv}; %\addlegendentry{AS};

        \addplot [domain=1E3:1E4] {1E-4/pow(x/1E4,2/2)} node[above, yshift=-2pt, midway, anchor=south west] {$h^{2}$};
        \addplot [domain=1E2:1E3] {3E-6/pow(x/1E3,3/2)} node[above, yshift=-2pt, midway, anchor=south west] {$h^{3}$};
      \end{loglogaxis}
    \end{tikzpicture}
    \caption{$\|\mathbf{u} - \mathbf{u}_h\|_{L^2}$}
  \end{subfigure}
  \begin{subfigure}[t]{0.45\textwidth} \centering
    \begin{tikzpicture}[scale=0.6]
       \begin{loglogaxis}[xlabel=Num DOFs, legend
           style={font=\LARGE},legend pos = outer north east,
          tick label style={font=\Large},
          label style={font=\Large},           
          ]
        \addplot table[x=dofs, y=pressureL2, col sep=comma]{code/stokes/data/CG2xCG1.2d.csv}; \addlegendentry{TH};
        \addplot table[x=dofs, y=pressureL2, col sep=comma]{code/stokes/data/CG1isoxCG1.2d.csv}; \addlegendentry{ISO};
        \addplot table[x=dofs, y=pressureL2, col sep=comma]{code/stokes/data/CG2xDG1alfeld.2d.csv}; \addlegendentry{SV};
        \addplot table[x=dofs, y=pressureL2, col sep=comma]{code/stokes/data/JM1xDG1.2d.csv}; \addlegendentry{JM};
        \addplot table[x=dofs, y=pressureL2, col sep=comma]{code/stokes/data/GN1xDG0.2d.csv}; \addlegendentry{GN};
        \addplot table[x=dofs, y=pressureL2, col sep=comma]{code/stokes/data/AS2xCG1.2d.csv}; \addlegendentry{AS};

        \addplot [domain=1E3:1E4] {2E-2/pow(x/1E4,1/2)} node[above, yshift=-2pt, midway, anchor=south west] {$h^{1}$};
        \addplot [domain=1E2:1E3] {6E-4/pow(x/1E3,2/2)} node[above, yshift=-2pt, midway, anchor=south west] {$h^{2}$};
      \end{loglogaxis}
    \end{tikzpicture}
    \caption{$\|p - p_h\|_{L^2}$}
  \end{subfigure}
  
  \begin{subfigure}[t]{0.45\textwidth} \centering
   \begin{tikzpicture}[scale=0.6]
     \begin{loglogaxis}[xlabel=Num DOFs,
         legend style={font=\LARGE},legend pos = outer north east,
         tick label style={font=\Large},
         label style={font=\Large},         
        ]
        \addplot table[x=dofs, y=stressL2, col sep=comma]{code/stokes/data/CG2xCG1.2d.csv}; %\addlegendentry{TH};
        \addplot table[x=dofs, y=stressL2, col sep=comma]{code/stokes/data/CG1isoxCG1.2d.csv}; %\addlegendentry{ISO};
        \addplot table[x=dofs, y=stressL2, col sep=comma]{code/stokes/data/CG2xDG1alfeld.2d.csv}; %\addlegendentry{SV};
        \addplot table[x=dofs, y=stressL2, col sep=comma]{code/stokes/data/JM1xDG1.2d.csv}; %\addlegendentry{JM};
        \addplot table[x=dofs, y=stressL2, col sep=comma]{code/stokes/data/GN1xDG0.2d.csv}; %\addlegendentry{GN};
        \addplot table[x=dofs, y=stressL2, col sep=comma]{code/stokes/data/AS2xCG1.2d.csv}; %\addlegendentry{AS};

        \addplot [domain=1E3:1E4] {3E-2/pow(x/1E4,1/2)} node[above, yshift=-2pt, midway, anchor=south west] {$h^{1}$};
        \addplot [domain=1E2:1E3] {1E-3/pow(x/1E3,2/2)} node[above, yshift=-2pt, midway, anchor=south west] {$h^{2}$};
      \end{loglogaxis}
    \end{tikzpicture}
    \caption{$\|\sigma - \sigma_h\|_{L^2}$}
  \end{subfigure}  
  \begin{subfigure}[t]{0.45\textwidth} \centering
    \begin{tikzpicture}[scale=0.6]
      \begin{loglogaxis}[xlabel=Num DOFs, legend
          style={font=\LARGE},legend pos = outer north east,
          tick label style={font=\Large},
          label style={font=\Large},          
        ]
        \addplot table[x=dofs, y=divL2, col sep=comma]{code/stokes/data/CG2xCG1.2d.csv}; \addlegendentry{TH};
        \addplot table[x=dofs, y=divL2, col sep=comma]{code/stokes/data/CG1isoxCG1.2d.csv}; \addlegendentry{ISO};
        \addplot table[x=dofs, y=divL2, col sep=comma]{code/stokes/data/CG2xDG1alfeld.2d.csv}; \addlegendentry{SV};
        \addplot table[x=dofs, y=divL2, col sep=comma]{code/stokes/data/JM1xDG1.2d.csv}; \addlegendentry{JM};
        \addplot table[x=dofs, y=divL2, col sep=comma]{code/stokes/data/GN1xDG0.2d.csv}; \addlegendentry{GN};
        \addplot table[x=dofs, y=divL2, col sep=comma]{code/stokes/data/AS2xCG1.2d.csv}; \addlegendentry{AS};

        \addplot [domain=1E2:1E4] {5E-2/pow(x/625,1/2)} node[above, yshift=-2pt, midway, anchor=south west] {$h^{1}$};
        \addplot [domain=1E2:1E4] {4E-5/pow(x/625,2/2)} node[above, yshift=-16pt, midway, anchor=south west] {$h^{2}$};
      \end{loglogaxis}
    \end{tikzpicture}
     \caption{$\|\div \mathbf{u}_h\|_{L^2}$}
  \end{subfigure}   

  \caption{Convergence plots for various Stokes discretizations on an $N \times N$ mesh of squares divided into right triangles.
  Taylor--Hood, Scott--Vogelius, Guzm\'an--Neilan, and Alfeld--Sorokina use the velocity-pressure formulation~\eqref{eq:velprestokes},
  while Johnson--Mercier uses the stress-velocity formulation~\eqref{eq:stressvelstokes}.
  Taylor--Hood elements are used on the original mesh, while the other formulations use macro-elements based on the Alfeld split.}
  \Description{Convergence plots for triangular Stokes discretizations.  The $L^2$ error in pressure, velocity, and divergence of velocity are plotting against the total number of degrees of freedom.}
  \label{fig:stokesconv2dvsverts}
\end{figure}

\Cref{fig:stokesconv2dvsverts} shows the error in the velocity, pressure, stress, and divergence versus the number of degrees of freedom  on each mesh.
In the velocity/pressure formulation, we solve the system and then compute a stress approximation via $\sigma = \epsilon(u) - p I$, and we similarly postprocess the stress to find the pressure in the stress-velocity formulation.
For both formulations we set $\nu=1$, and the theoretically predicted rates of convergence are observed.
Taylor--Hood method actually gives the best velocity and stress approximations on a given mesh, followed by the Alfeld-Sorokina method, with these two also giving the best approximation of the pressure.
The Guzm\'an--Neilan approach uses a lower-order pressure approximation, which limits the overall accuracy of the method.
The Taylor--Hood, iso, and Johnson--Mercier methods have quite a large residual divergence, while the other three approaches enforce the divergence-free condition to high accuracy.  
These results indicate that Taylor--Hood may be preferred if a pointwise divergence-free condition is not critical, while Alfeld--Sorokina may be preferred if it is.
The stress--velocity formulation with Johnson--Mercier elements is also competitive and may be of interest for more complex rheology.

\Cref{fig:stokesconv3dvsverts} shows analogous results for three dimensions.

The minimal Scott--Vogelius method is higher order, with cubic velocities and discontinuous quadratic pressures on the Alfeld split, giving much lower error at the greater expense.  Otherwise, as in 2D, Taylor--Hood and Alfeld--Sorokina give the most accurate answers, differing most in the divergence-free condition.

\begin{figure}[htbp]
  \begin{subfigure}[t]{0.45\textwidth} \centering
    \begin{tikzpicture}[scale=0.6]
      \begin{loglogaxis}[xlabel=Num DOFs, legend style={font=\LARGE},
          legend pos = outer north east,
          tick label style={font=\Large},
          label style={font=\Large},          
         ]
        \addplot table[x=dofs, y=velocityL2, col sep=comma]{code/stokes/data/CG2xCG1.3d.csv};            %\addlegendentry{TH};
        \addplot table[x=dofs, y=velocityL2, col sep=comma]{code/stokes/data/CG1isoxCG1.3d.csv};         %\addlegendentry{ISO};
        \addplot table[x=dofs, y=velocityL2, col sep=comma]{code/stokes/data/CG3xDG2alfeld.3d.csv};      %\addlegendentry{SV};
        \addplot table[x=dofs, y=velocityL2, col sep=comma]{code/stokes/data/JM1xDG1.3d.csv};            %\addlegendentry{JM};
        \addplot table[x=dofs, y=velocityL2, col sep=comma]{code/stokes/data/GN1xDG0.3d.csv};            %\addlegendentry{GN};
        \addplot table[x=dofs, y=velocityL2, col sep=comma]{code/stokes/data/GNH1div3xCG1alfeld.3d.csv}; %\addlegendentry{AS};

        \addplot [domain=1E4:1E6] {5E-3/pow(x/1E4,2/3)} node[above, yshift=-2pt, midway, anchor=south west] {$h^{2}$};
        \addplot [domain=1E4:1E6] {1E-5/pow(x/1E5,3/3)} node[above, yshift=-2pt, midway, anchor=south west] {$h^{3}$};
        \addplot [domain=1E5:1E7] {5E-6/pow(x/1E5,4/3)} node[above, yshift=-2pt, midway, anchor=south west] {$h^{4}$};
      \end{loglogaxis}
    \end{tikzpicture}
    \caption{$\|\mathbf{u} - \mathbf{u}_h\|_{L^2}$}
  \end{subfigure}
  \begin{subfigure}[t]{0.45\textwidth} \centering
    \begin{tikzpicture}[scale=0.6]
       \begin{loglogaxis}[xlabel=Num DOFs, legend style={font=\LARGE},
          legend pos = outer north east,
          tick label style={font=\Large},
          label style={font=\Large},           
          ]
        \addplot table[x=dofs, y=pressureL2, col sep=comma]{code/stokes/data/CG2xCG1.3d.csv};            \addlegendentry{TH};
        \addplot table[x=dofs, y=pressureL2, col sep=comma]{code/stokes/data/CG1isoxCG1.3d.csv};         \addlegendentry{ISO};
        \addplot table[x=dofs, y=pressureL2, col sep=comma]{code/stokes/data/CG3xDG2alfeld.3d.csv};      \addlegendentry{SV};
        \addplot table[x=dofs, y=pressureL2, col sep=comma]{code/stokes/data/JM1xDG1.3d.csv};            \addlegendentry{JM};
        \addplot table[x=dofs, y=pressureL2, col sep=comma]{code/stokes/data/GN1xDG0.3d.csv};            \addlegendentry{GN};
        \addplot table[x=dofs, y=pressureL2, col sep=comma]{code/stokes/data/GNH1div3xCG1alfeld.3d.csv}; \addlegendentry{AS};

        \addplot [domain=1E4:1E6] {4E-2/pow(x/3E5,1/3)} node[above, yshift=-2pt, midway, anchor=south west] {$h^{1}$};
        \addplot [domain=1E5:1E7] {1.5E-3/pow(x/3E6,2/3)} node[above, yshift=-2pt, midway, anchor=south west] {$h^{2}$};
        \addplot [domain=1E5:1E7] {1.5E-5/pow(x/1E7,3/3)} node[above, yshift=-2pt, midway, anchor=south west] {$h^{3}$};
      \end{loglogaxis}
    \end{tikzpicture}
    \caption{$\|p - p_h\|_{L^2}$}
  \end{subfigure}

  \begin{subfigure}[t]{0.45\textwidth} \centering
    \begin{tikzpicture}[scale=0.6]
      \begin{loglogaxis}[xlabel=Num DOFs, 
          legend style={font=\LARGE}
          ,legend pos = outer north east,
          tick label style={font=\Large},
          label style={font=\Large},          
         ]
        \addplot table[x=dofs, y=stressL2, col sep=comma]{code/stokes/data/CG2xCG1.3d.csv};            %\addlegendentry{TH};
        \addplot table[x=dofs, y=stressL2, col sep=comma]{code/stokes/data/CG1isoxCG1.3d.csv};         %\addlegendentry{ISO};
        \addplot table[x=dofs, y=stressL2, col sep=comma]{code/stokes/data/CG3xDG2alfeld.3d.csv};      %\addlegendentry{SV};
        \addplot table[x=dofs, y=stressL2, col sep=comma]{code/stokes/data/JM1xDG1.3d.csv};            %\addlegendentry{JM};
        \addplot table[x=dofs, y=stressL2, col sep=comma]{code/stokes/data/GN1xDG0.3d.csv};            %\addlegendentry{GN};
        \addplot table[x=dofs, y=stressL2, col sep=comma]{code/stokes/data/GNH1div3xCG1alfeld.3d.csv}; %\addlegendentry{AS};
        \addplot [domain=1E4:1E6] {7E-2/pow(x/3E5,1/3)} node[above, yshift=-2pt, midway, anchor=south west] {$h^{1}$};
        \addplot [domain=1E5:1E7] {3E-3/pow(x/3E6,2/3)} node[above, yshift=-2pt, midway, anchor=south west] {$h^{2}$};
        \addplot [domain=1E5:1E7] {2E-5/pow(x/1E7,3/3)} node[above, yshift=-2pt, midway, anchor=south west] {$h^{3}$};

      \end{loglogaxis}
    \end{tikzpicture}
    \caption{$\|\sigma - \sigma_h\|_{L^2}$}
  \end{subfigure}  
  \begin{subfigure}[t]{0.45\textwidth} \centering
    \begin{tikzpicture}[scale=0.6]
      \begin{loglogaxis}[xlabel=Num DOFs, legend
          style={font=\LARGE},legend pos = outer north east,
          tick label style={font=\Large},
          label style={font=\Large},          
         ]
        \addplot table[x=dofs, y=divL2, col sep=comma]{code/stokes/data/CG2xCG1.3d.csv};            \addlegendentry{TH};
        \addplot table[x=dofs, y=divL2, col sep=comma]{code/stokes/data/CG1isoxCG1.3d.csv};         \addlegendentry{ISO};
        \addplot table[x=dofs, y=divL2, col sep=comma]{code/stokes/data/CG3xDG2alfeld.3d.csv};      \addlegendentry{SV};
        \addplot table[x=dofs, y=divL2, col sep=comma]{code/stokes/data/JM1xDG1.3d.csv};            \addlegendentry{JM};
        \addplot table[x=dofs, y=divL2, col sep=comma]{code/stokes/data/GN1xDG0.3d.csv};            \addlegendentry{GN};
        \addplot table[x=dofs, y=divL2, col sep=comma]{code/stokes/data/GNH1div3xCG1alfeld.3d.csv}; \addlegendentry{AS};
        \addplot [domain=1E4:1E6] {4E-2/pow(x/25^3,1/3)} node[above, yshift=-2pt, midway, anchor=south west] {$h^{1}$};
        \addplot [domain=1E3:1E6] {3E-5/pow(x/25^3,2/3)} node[above, yshift=-16pt, midway, anchor=south west] {$h^{2}$};
      \end{loglogaxis}
    \end{tikzpicture}
    \caption{$\|\div\mathbf{u}_h\|_{L^2}$}
  \end{subfigure}   

  \caption{Convergence plots for various Stokes discretizations on an $N \times N \times N$ mesh of cubes divided into six tetrahedra.
  Taylor--Hood, Scott--Vogelius, Guzm\'an--Neilan, and Alfeld--Sorokina use the velocity-pressure formulation~\eqref{eq:velprestokes},
  while Johnson--Mercier uses the stress-velocity formulation~\eqref{eq:stressvelstokes}.  Taylor--Hood elements are used on the original mesh, while the other formulations use macro-elements based on the Alfeld split.}
  \Description{Convergence plots for tetrahedral Stokes discretizations.  The $L^2$ error in pressure, velocity, and divergence of velocity are plotting against the total number of degrees of freedom.}  
  \label{fig:stokesconv3dvsverts}
\end{figure}

\begin{figure}[htbp]
  \begin{subfigure}{0.49\textwidth}
    \begin{tikzpicture}[scale=0.8]
      \pgfplotstableread[col sep=comma,]{code/stokes/data/stokesflops2d.csv}\datatable
      \begin{axis}[ybar stacked,
          ylabel near ticks,
          ymin=1e3, ymax=0.4e5,
          % ymode = log,
        legend pos=north west,
        %% ymin=-1e4, ymax=0.5e4,
        legend cell align=left,
        legend style={draw=none},
        xtick=data,
        xticklabels from table={\datatable}{Name}, ylabel={FLOPs},
        xticklabel style={xshift=-10pt, rotate=60}]
        \addplot+ table [x expr=\coordindex, y={A}]{\datatable};
        \addplot+[pattern=north east lines, pattern color=red]
          table [x expr=\coordindex, y={B}]{\datatable};
        \legend{$A FLOPs$, $B FLOPs$}
      \end{axis}
    \end{tikzpicture}
    \caption{2D}
    \label{fig:stokesflops2d}
  \end{subfigure}
  \begin{subfigure}{0.49\textwidth}
    \begin{tikzpicture}[scale=0.8]
      \pgfplotstableread[col sep=comma,]{code/stokes/data/stokesflops3d.csv}\datatable
      \begin{axis}[ybar stacked,
          ylabel near ticks,
          %ymode = log,
          legend pos=north west,
          ymin=1e4, ymax=0.4e7,
        %% ymin=-1e4, ymax=0.5e4,
        legend cell align=left,
        legend style={draw=none},
        xtick=data,
        xticklabels from table={\datatable}{Name}, ylabel={FLOPs},
        xticklabel style={xshift=-10pt, rotate=60}]
        \addplot+ table [x expr=\coordindex, y={A}]{\datatable};
        \addplot+[pattern=north east lines, pattern color=red]
          table [x expr=\coordindex, y={B}]{\datatable};
        \legend{$A FLOPs$, $B FLOPs$}
      \end{axis}
    \end{tikzpicture}
    \caption{3D}
    \label{fig:stokesflops3d}
  \end{subfigure}
  
    \caption{FLOP count for evaluating the kernels for computing the
      elementwise contributions to matrices $A$ and $B$ from~\eqref{eq:stokesblockmat} for the Stokes operator.  This count assumes reference basis elements are pre-tabulated and then includes the cost of transforming basis functions and their derivatives and performing integration.}
    \Description{Bar chart showing the number of floating point operations required to form the element matrices for various Stokes elements.} 
    \label{fig:stokesflops}
\end{figure}
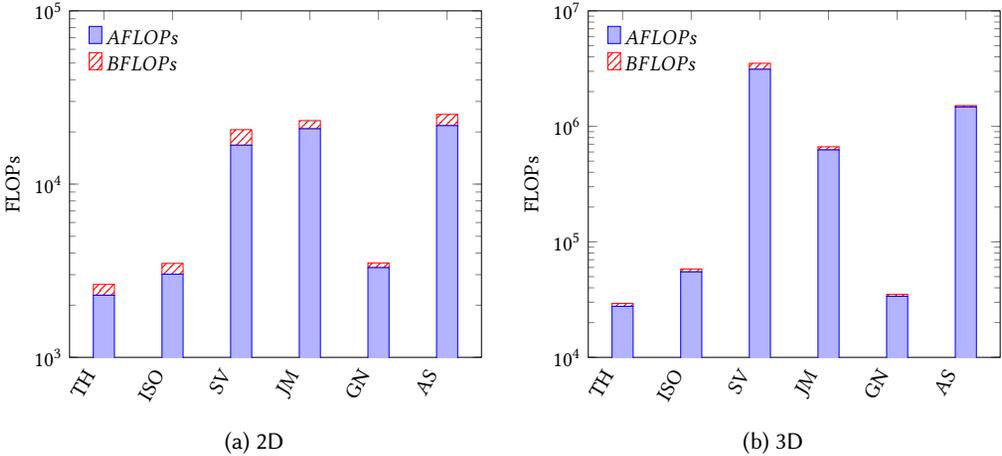

We also report the number of floating-point operations (FLOPs) required to assemble the element-level matrices.  The formulations~\eqref{eq:velprestokes} and~\eqref{eq:stressvelstokes} both lead to block matrices of the form
\begin{equation}
  \label{eq:stokesblockmat}
  \begin{bmatrix}
    A & B^\top \\
    B & 0 
  \end{bmatrix},
\end{equation}
where the matrices $A$ and $B$ differ between the pressure-velocity and stress-velocity formulations.  In either case, both matrices are formed by iterating over cells to form local contributions and assembling them into the global sparse matrix.  \Cref{fig:stokesflops} gives the operation count reported by \lstinline{tsfc} in the element-level kernels for $A$ and $B$.
These numbers assume that the reference basis functions are pre-tabulated but include the entire cost of forming Jacobians, transforming basis functions, and integrating over the cell.
In each case, the cost of forming $A$ greatly dominates that of $B$.  We see that the non-macro nature of the Taylor--Hood pair makes assembly much cheaper than many of the alternatives.
Note that, although the ISO element pair uses a macroelement velocity, it only uses piecewise linear polynomials on the splitting, and hence $A$ requires only a single quadrature point on each subcell.
Since total run-time typically depends more strongly on algebraic solvers than assembly, the advantages of macroelements may outweigh the increased assembly costs.

Per the discussion above, the FLOP count to form element-level matrices can vary greatly among our different macroelements.
For the Stokes equations, having a macroelement pressure space tends to increase the dimension of that space and also make each velocity basis function couple to more pressure basis functions, for example.
Table~\ref{table:stokessparsity} shows basic sparsity information for several Stokes pairs. 
The sparsity structure of the assembled matrix depends somewhat on macroelement structure, but also on the global degrees of freedom and their location -- vertex degrees of freedom are shared by more cells and so they tend to reduce the global number of degrees of freedom while creating more couplings (that is, reducing sparsity).

\begin{table}
  \begin{tabular}{ccc} 
    Element & Rows & Nonzeros per row \\ \hline
    TH & 659 & 18.64 \\
    ISO & 659 & 18.64 \\
    SV & 2754 & 19.49 \\
    JM & 1984 & 21.73 \\
    GN & 498 & 16.04 \\
    AS & 868 & 29.70 \\
  \end{tabular}
  \caption{The number of degrees of freedom (rows in the matrix) and average sparsity for the Stokes equations discretized on an $8\times 8$ mesh divided into triangles.}
  \label{table:stokessparsity}
\end{table}

\begin{figure}[htbp]
  \phantom{.}\hfill
  \begin{subfigure}[t]{0.15\textwidth}
    \centering
    \includegraphics[height=0.95\textwidth]{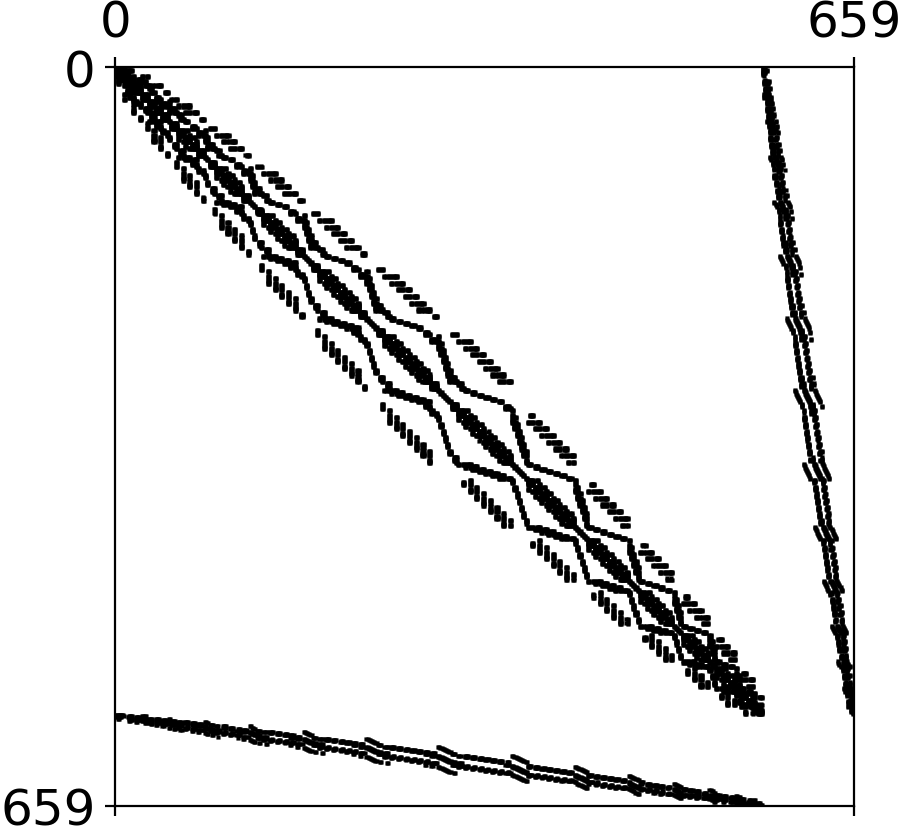}
    \caption{TH}
  \end{subfigure}
   \hfill
  \begin{subfigure}[t]{0.15\textwidth}
    \centering
    \includegraphics[height=0.95\textwidth]{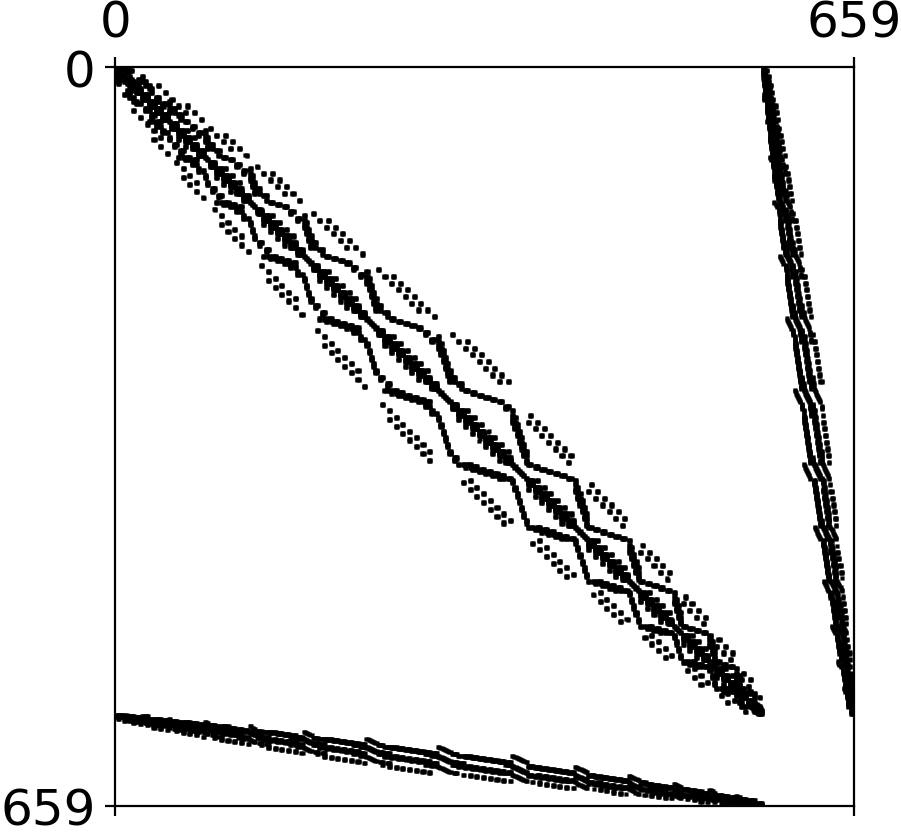}
    \caption{ISO}
  \end{subfigure}
   \hfill
  \begin{subfigure}[t]{0.15\textwidth}
    \centering
    \includegraphics[height=0.95\textwidth]{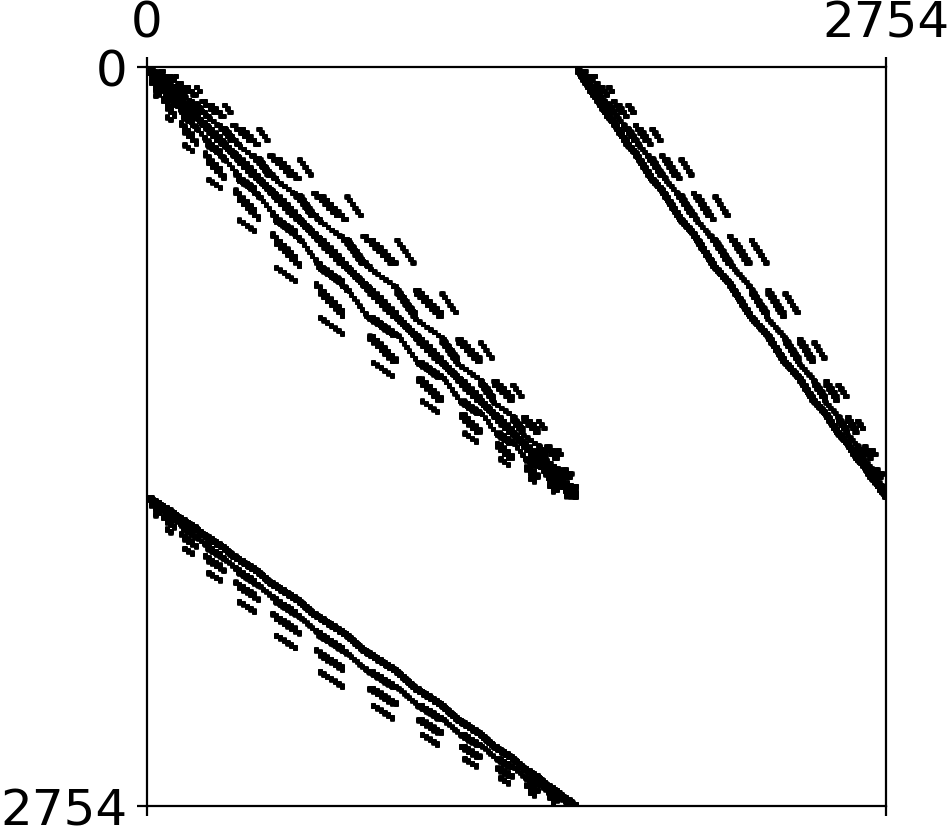}
    \caption{SV}
  \end{subfigure}
   \hfill
  \begin{subfigure}[t]{0.15\textwidth}
    \centering
    \includegraphics[height=0.95\textwidth]{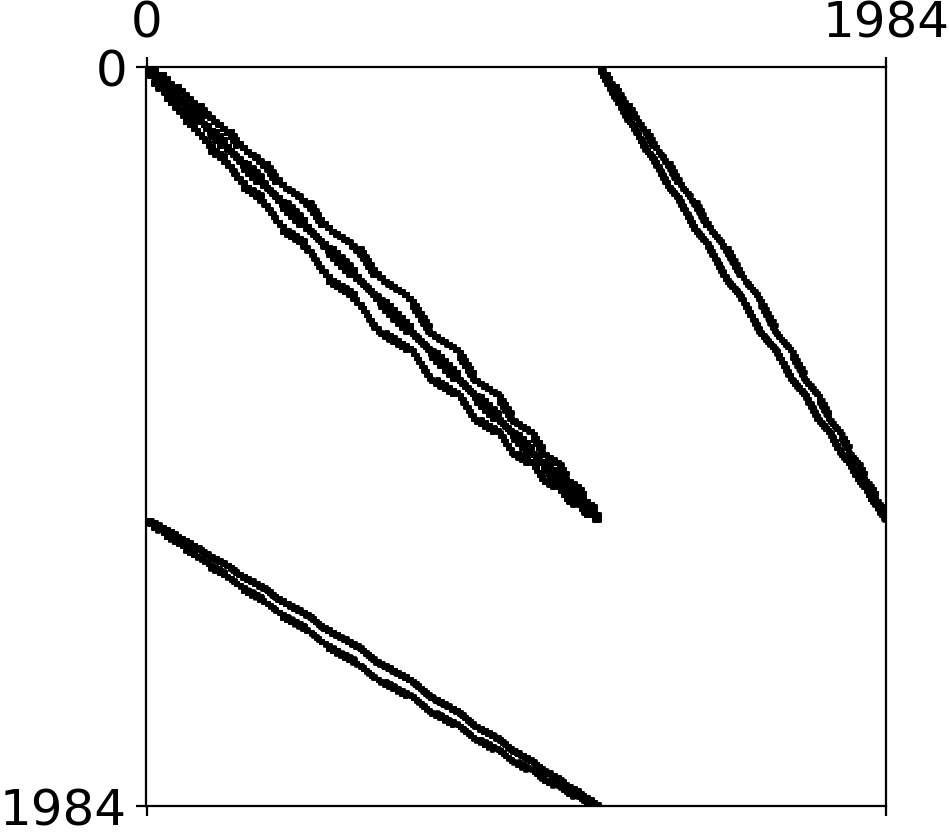}
    \caption{JM}
  \end{subfigure}
   \hfill
  \begin{subfigure}[t]{0.15\textwidth}
    \centering
    \includegraphics[height=0.95\textwidth]{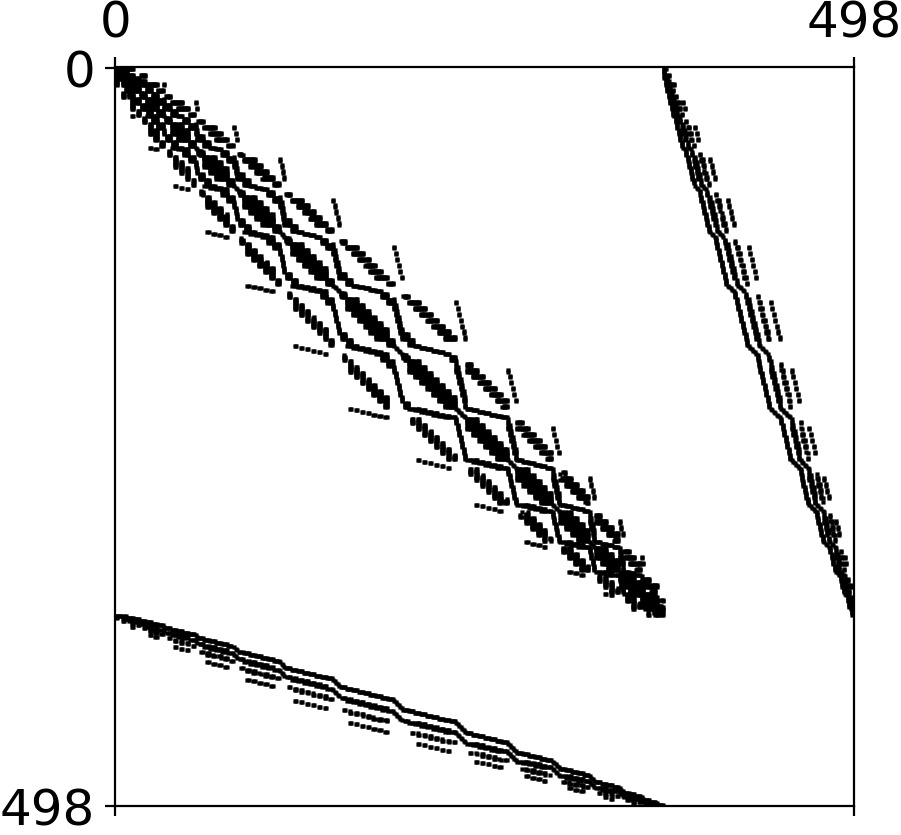}
    \caption{GN}
  \end{subfigure}
   \hfill
  \begin{subfigure}[t]{0.15\textwidth}
    \centering
    \includegraphics[height=0.95\textwidth]{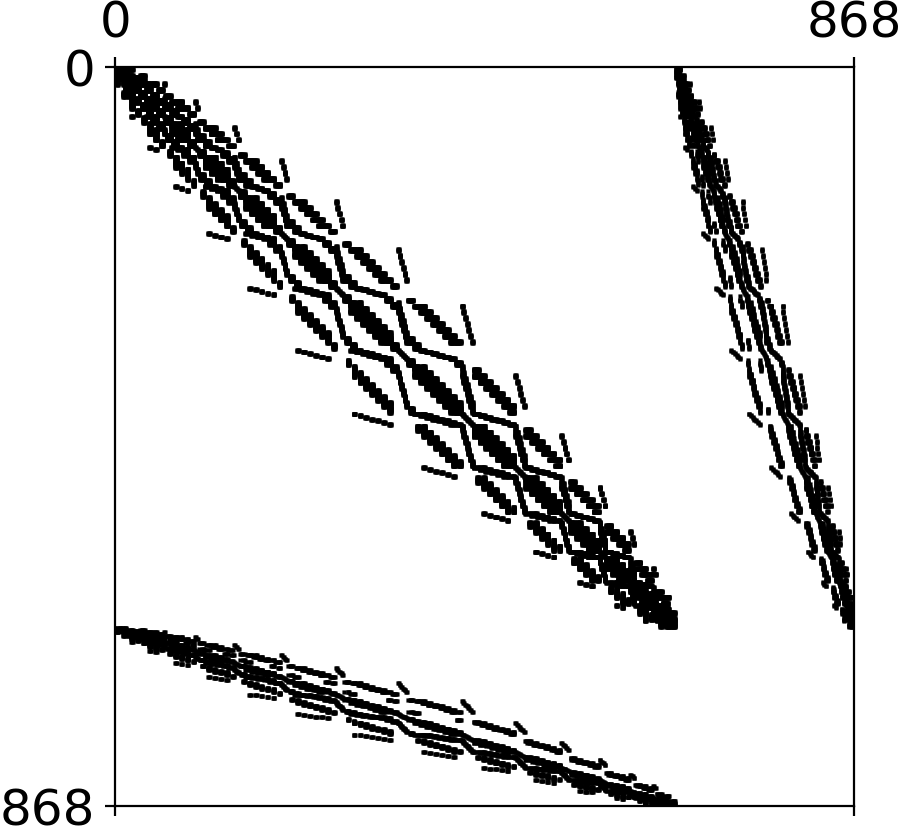}
    \caption{AS}
  \end{subfigure}  
  \hfill\phantom{.}
  \caption{Sparsity patterns for several Stokes pairs on an $8 \times 8$ mesh subdivided into right triangles.}
  \Description{Plots of the sparsity patterns for the Stokes operator on an $8 \times 8$ mesh using various discretizations.}  
  \label{fig:stokessparsity}
\end{figure}

We also report on the relative time for solving and assembling linear systems on a range of meshes used in our convergence study above.  A full treatment of solvers is beyond the scope of this paper, but we have chosen a relatively simple configuration -- we simply build the matrix and use the MUMPS LU factorization.
Figure~\ref{fig:stokesruntime} shows a slightly superlinear timing versus the number of degrees of freedom, and Figure~\ref{fig:stokesassemblyfraction} shows that the total assembly time does not exceed about 20\% of the run-time.
For problems at the scale we consider, sparse direct is likely to be quite efficient compared to iterative methods, although we have not performed an exhaustive comparison.

\begin{figure}[htbp]
  \begin{subfigure}[t]{0.45\textwidth} \centering
    \begin{tikzpicture}[scale=0.6]
      \begin{loglogaxis}[xlabel=Num DOFs, legend style={font=\LARGE},
          tick label style={font=\Large},
          label style={font=\Large}, ylabel={Time (s)}
         ]
       \addplot table[x=dofs, y=SNESSolve, col sep=comma]{code/stokes/data/CG2xCG1.2d.csv};
        \addplot table[x=dofs, y=SNESSolve, col sep=comma]{code/stokes/data/CG1isoxCG1.2d.csv};
        \addplot table[x=dofs, y=SNESSolve, col sep=comma]{code/stokes/data/CG2xDG1alfeld.2d.csv};
        \addplot table[x=dofs, y=SNESSolve, col sep=comma]{code/stokes/data/JM1xDG1.2d.csv};
        \addplot table[x=dofs, y=SNESSolve, col sep=comma]{code/stokes/data/GN1xDG0.2d.csv};
        \addplot table[x=dofs, y=SNESSolve, col sep=comma]{code/stokes/data/AS2xCG1alfeld.2d.csv};
        \addplot [domain=4E3:3E4] {1E-2 * pow(x/4E3,1)} node[above, xshift=16pt, yshift=-8pt, midway] {$\mathcal{O}(N)$};        
      \end{loglogaxis}
    \end{tikzpicture}
    \caption{Total run-time}
    \label{fig:stokesruntime}
  \end{subfigure}
  \begin{subfigure}[t]{0.45\textwidth} \centering
    \begin{tikzpicture}[scale=0.6]
      \begin{semilogxaxis}[xlabel=Num DOFs, legend
          style={font=\LARGE},legend pos = outer north east,
          tick label style={font=\Large},
          label style={font=\Large}, ymin=0, ymax=1,
          ylabel={Assembly fraction}
        ]
        \addplot table[x=dofs, y expr= \thisrow{SNESJacobianEval} / \thisrow{SNESSolve}, col sep=comma]{code/stokes/data/CG2xCG1.2d.csv}; \addlegendentry{TH};
        \addplot table[x=dofs, y expr= \thisrow{SNESJacobianEval} / \thisrow{SNESSolve}, col sep=comma]{code/stokes/data/CG1isoxCG1.2d.csv}; \addlegendentry{ISO}; 
        \addplot table[x=dofs, y expr= \thisrow{SNESJacobianEval} / \thisrow{SNESSolve}, col sep=comma]{code/stokes/data/CG2xDG1alfeld.2d.csv}; \addlegendentry{SV};
        \addplot table[x=dofs, y expr= \thisrow{SNESJacobianEval} / \thisrow{SNESSolve}, col sep=comma]{code/stokes/data/JM1xDG1.2d.csv}; \addlegendentry{JM};
        \addplot table[x=dofs, y expr= \thisrow{SNESJacobianEval} / \thisrow{SNESSolve}, col sep=comma]{code/stokes/data/GN1xDG0.2d.csv}; \addlegendentry{GN};
        \addplot table[x=dofs, y expr= \thisrow{SNESJacobianEval} / \thisrow{SNESSolve}, col sep=comma]{code/stokes/data/AS2xCG1alfeld.2d.csv}; \addlegendentry{AS};        
    \end{semilogxaxis}
    \end{tikzpicture}
    \caption{Fraction of time in assembly}
    \label{fig:stokesassemblyfraction}
  \end{subfigure}
  \caption{Timing results versus the number of degrees of freedom for the two-dimensional Stokes equations with various element pairs.  We show the total solve time (to assemble, factor, and solve) each system on the left and the fraction of that time spent in assembly on the right.}
  \Description{Plots showing the run-time and ratio of time spent in assembly to sparse direct factorization in some two-dimensional Stokes calculations.}
  \label{fig:timingstokes}
\end{figure}

\subsection{Navier-Stokes}
We also applied our methods to a well-known benchmark for the two-dimensional steady-state Navier-Stokes equations, measuring the drag and lift on a cylinder and the pressure drop across it~\cite{john2004reference,schafer1996benchmark}.
The domain is given by $\Omega=[0,2.2] \times [0, 0.41] \backslash B_r(0.2, 0.2)$, with radius $r=0.05$, and is shown in \Cref{fig:domain2d}. The density is taken as $\rho=1$ and kinematic viscosity is $\nu = 10^{-3}$, which gives a Reynolds number of 20.
No-slip conditions are imposed on the top and bottom of the pipe and the cylinder.  Natural boundary conditions (no-stress) are imposed on the outflow right end, and a parabolic profile is posed on the inflow boundary on the left end:
\begin{equation}
u(0, y) = \left( \frac{4 y(0.41 - y)}{0.41^2}, 0 \right) \equiv U_\mathrm{in}(y).
\end{equation}

\begin{figure}[htbp]
\begin{subfigure}[t]{0.45\textwidth}
  \begin{tikzpicture}[scale=0.6]
      \begin{loglogaxis}[xlabel=Num DOFs, legend style={font=\LARGE},
          legend pos = outer north east,
          tick label style={font=\Large},
          label style={font=\Large}
         ]
         \addplot table[x=dofs, y expr={abs(\thisrow{drag} - 5.57953523384)}, col sep=comma] {code/nse/data/CG2xCG1.2d.csv};         %\addlegendentry{TH};
         \addplot table[x=dofs, y expr={abs(\thisrow{drag} - 5.57953523384)}, col sep=comma] {code/nse/data/CG1isoxCG1.2d.csv};      %\addlegendentry{ISO};
         \addplot table[x=dofs, y expr={abs(\thisrow{drag} - 5.57953523384)}, col sep=comma] {code/nse/data/CG2xDG1alfeld.2d.csv};   %\addlegendentry{SV};        
         \addplot table[x=dofs, y expr={abs(\thisrow{drag} - 5.57953523384)}, col sep=comma] {code/nse/data/JM1xDG1.2d.csv};   %\addlegendentry{JM};        
         \addplot table[x=dofs, y expr={abs(\thisrow{drag} - 5.57953523384)}, col sep=comma] {code/nse/data/GN1xDG0.2d.csv};         %\addlegendentry{GN};
         \addplot table[x=dofs, y expr={abs(\thisrow{drag} - 5.57953523384)}, col sep=comma] {code/nse/data/AS2xCG1.2d.csv};         %\addlegendentry{AS};        

        \addplot [domain=6.8e4:1.1e6] {8E-1/pow(x/1E5,1/2)} node[above, yshift=-2pt, midway, anchor=south west] {$h^{1}$};
        \addplot [domain=6.8e3:1.1e5] {1E-3/pow(x/1E5,2/2)} node[above, yshift=-2pt, midway, anchor=south west] {$h^{2}$};
      \end{loglogaxis}
  \end{tikzpicture}
  \caption{Drag error}
\end{subfigure}
\begin{subfigure}[t]{0.45\textwidth}
  \begin{tikzpicture}[scale=0.6]
      \begin{loglogaxis}[xlabel=Num DOFs, legend style={font=\LARGE},
          legend pos = outer north east,
          tick label style={font=\Large},
          label style={font=\Large}          
         ]
         \addplot table[x=dofs, y expr={abs(\thisrow{lift} - 0.010618948146)}, col sep=comma] {code/nse/data/CG2xCG1.2d.csv};       \addlegendentry{TH};
         \addplot table[x=dofs, y expr={abs(\thisrow{lift} - 0.010618948146)}, col sep=comma] {code/nse/data/CG1isoxCG1.2d.csv};    \addlegendentry{ISO};
         \addplot table[x=dofs, y expr={abs(\thisrow{lift} - 0.010618948146)}, col sep=comma] {code/nse/data/CG2xDG1alfeld.2d.csv}; \addlegendentry{SV};        
         \addplot table[x=dofs, y expr={abs(\thisrow{lift} - 0.010618948146)}, col sep=comma] {code/nse/data/JM1xDG1.2d.csv}; \addlegendentry{JM};        
         \addplot table[x=dofs, y expr={abs(\thisrow{lift} - 0.010618948146)}, col sep=comma] {code/nse/data/GN1xDG0.2d.csv};       \addlegendentry{GN};        
         \addplot table[x=dofs, y expr={abs(\thisrow{lift} - 0.010618948146)}, col sep=comma] {code/nse/data/AS2xCG1.2d.csv};       \addlegendentry{AS};        
        
        \addplot [domain=6.8e4:1.1e6] {4E-2/pow(x/1E5,1/2)} node[above, yshift=-2pt, midway, anchor=south west] {$h^{1}$};
        \addplot [domain=6.8e3:1.1e5] {2E-5/pow(x/1E5,2/2)} node[above, yshift=-22pt, midway, anchor=south west] {$h^{2}$};
      \end{loglogaxis}
  \end{tikzpicture}
  \caption{Lift error}
\end{subfigure}

\begin{subfigure}[t]{0.45\textwidth}
  \begin{tikzpicture}[scale=0.6]
      \begin{loglogaxis}[xlabel=Num DOFs, legend style={font=\LARGE},
          legend pos = outer north east,
          tick label style={font=\Large},
          label style={font=\Large}          
         ]
         \addplot table[x=dofs, y expr={abs(\thisrow{pressureDrop} - 0.11752016697)}, col sep=comma] {code/nse/data/CG2xCG1.2d.csv};       %\addlegendentry{TH};
         \addplot table[x=dofs, y expr={abs(\thisrow{pressureDrop} - 0.11752016697)}, col sep=comma] {code/nse/data/CG1isoxCG1.2d.csv};    %\addlegendentry{ISO};
         \addplot table[x=dofs, y expr={abs(\thisrow{pressureDrop} - 0.11752016697)}, col sep=comma] {code/nse/data/CG2xDG1alfeld.2d.csv}; %\addlegendentry{SV};  
         \addplot table[x=dofs, y expr={abs(\thisrow{pressureDrop} - 0.11752016697)}, col sep=comma] {code/nse/data/JM1xDG1.2d.csv}; %\addlegendentry{JM};       
         \addplot table[x=dofs, y expr={abs(\thisrow{pressureDrop} - 0.11752016697)}, col sep=comma] {code/nse/data/GN1xDG0.2d.csv};       %\addlegendentry{GN};       
         \addplot table[x=dofs, y expr={abs(\thisrow{pressureDrop} - 0.11752016697)}, col sep=comma] {code/nse/data/AS2xCG1.2d.csv};       %\addlegendentry{AS};       

        \addplot [domain=6.1e4:1.1e6] {7E-3/pow(x/1E5,1/2)} node[above, yshift=-2pt, midway, anchor=south west] {$h^{1}$};
        \addplot [domain=6.1e3:1.1e5] {1E-5/pow(x/1E5,2/2)} node[above, yshift=-2pt, midway, anchor=south west] {$h^{2}$};
      \end{loglogaxis}
  \end{tikzpicture}
  \caption{Pressure drop error}
\end{subfigure}
   \begin{subfigure}[t]{0.45\textwidth}
  \begin{tikzpicture}[scale=0.6]
      \begin{loglogaxis}[xlabel=Num DOFs, legend style={font=\LARGE},
          legend pos = outer north east,
          tick label style={font=\Large},
          label style={font=\Large}          
         ]
         \addplot table[x=dofs, y=errorDiv, col sep=comma] {code/nse/data/CG2xCG1.2d.csv};          \addlegendentry{TH};
         \addplot table[x=dofs, y=errorDiv, col sep=comma] {code/nse/data/CG1isoxCG1.2d.csv};       \addlegendentry{ISO};
         \addplot table[x=dofs, y=errorDiv , col sep=comma] {code/nse/data/CG2xDG1alfeld.2d.csv};   \addlegendentry{SV};        
         \addplot table[x=dofs, y=errorDiv, col sep=comma] {code/nse/data/JM1xDG1.2d.csv};    \addlegendentry{JM};        
         \addplot table[x=dofs, y=errorDiv , col sep=comma] {code/nse/data/GN1xDG0.2d.csv};         \addlegendentry{GN};        
         \addplot table[x=dofs, y=errorDiv , col sep=comma] {code/nse/data/AS2xCG1.2d.csv};         \addlegendentry{AS};        
        
        \addplot [domain=6.8e4:1.1e6] {5E-1/pow(x/1E5,1/2)} node[above, yshift=-2pt, midway, anchor=south west] {$h^{1}$};
        \addplot [domain=6.8e3:1.1e5] {1E-3/pow(x/1E5,2/2)} node[above, yshift=-16pt, midway, anchor=south west] {$h^{2}$};
      \end{loglogaxis}
  \end{tikzpicture}
  \caption{$\|\div\mathbf{u}_h\|_{L^2}$}
  \label{fig:2dnsediv}
\end{subfigure}

\caption{Error in drag, lift, pressure drop, and divergence for the 2d cylinder problem.
   Here, the Johnson--Mercier discretization of the stress-velocity formulation
   outperforms the Taylor--Hood and Scott--Vogelius discretization of the
   velocity-pressure formulation on each refinement level.}
\Description{Convergence plots for the Navier--Stokes benchmark.  These show the error in drag, lift, pressure drop, and divergence of velocity versus the number of degrees of freedom.}
\label{fig:2ddraglift}
\end{figure}
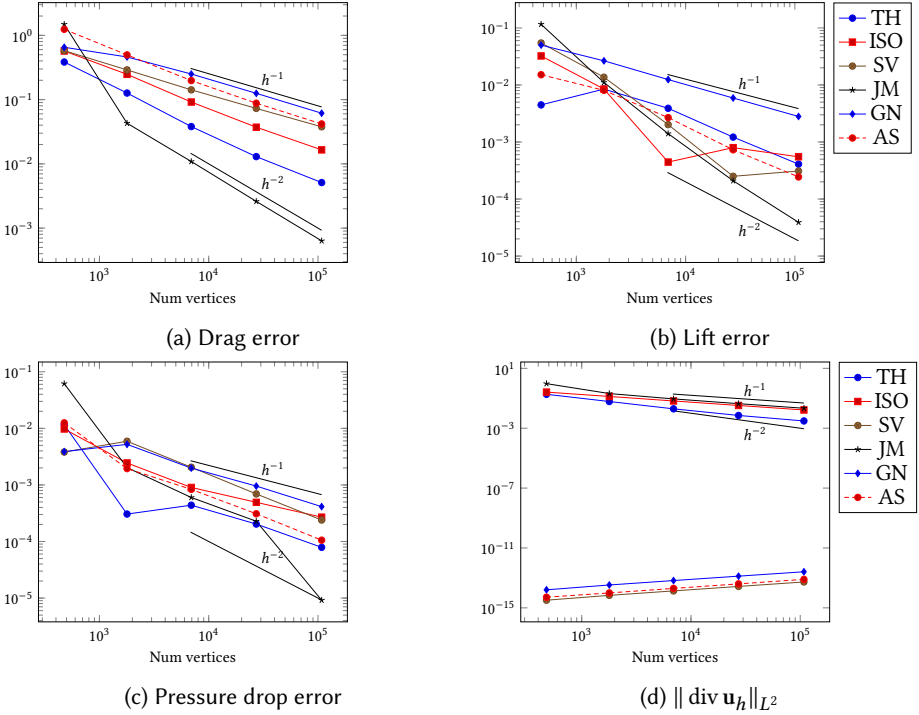

\Cref{fig:2ddraglift} shows the error in several key quantities for each of the methods considered for Stokes flow.
These quantities are the lift and drag on the cylinder and the difference in pressure $p(0.15, 0.2) - p(0.25, 0.2)$ across the cylinder.
These values are known to high accuracy~\cite{nabh1998high}.
We also report the norm of the divergence of the velocity per~\eqref{eq:divnorm}.

\begin{figure}[htbp]
  \centering
  \begin{tikzpicture}[scale=4]
    \draw (0,0) rectangle (2.2, 0.41);
    \draw (.2, .2) circle (0.05);
    \node[] at (1.1, 0.46) {$u=0$};
    \node[] at (1.1, -0.05) {$u=0$};
    \node[] at (.35, .2) {$u=0$};
    \node[] at (2.4, 0.205) {$\sigma n = 0$};
     \node[] at (-0.3, 0.205) {$u=U_\mathrm{in}(y)$};
  \end{tikzpicture}
  \caption{Computational domain for flow past cylinder, with boundary conditions indicated on each part of the boundary.}
  \Description{Schematic drawing of the computational domain -- a rectangular with a circle removed -- used in the benchmark calculations.}
  \label{fig:domain2d}
\end{figure}
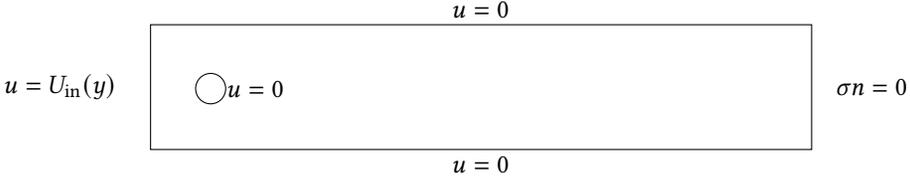

Perhaps surprisingly, we obtain a different relative ordering of accuracy between the various methods than for the Stokes equations, and this varies depending on the quantity of interest.
The Johnson--Mercier formulation gives the best drag approximation, and nearly the best results for lift and pressure drop.  The Taylor--Hood pair is competitive (and nearly the best for drag), but neither is it divergence free.
Among the divergence-free methods, we see that Scott--Vogelius and Alfeld--Sorokina pairs are more accurate than the lower-order Guzm\'an--Neilan pair.
It is conceivable that these relative orderings would change at higher Reynolds number (for example, Taylor--Hood is not pressure robust).
That each problem and functional may be best resolved by a different element (and even variational formulation) supports our goal of enabling a quite broad class of discretizations in the Firedrake code stack.

\subsection{Fourth-order problems}
Next, we apply our newly-enabled  $C^1$ macroelements to the plate-bending biharmonic problem,
\begin{equation}
  \Delta^2 u = f
\end{equation}
on some domain $\Omega \subset \mathbb{R}^2$.
In our examples, we consider clamped boundary
conditions $u = \tfrac{\partial u}{\partial n} = 0$ on $\partial \Omega$.
The obvious bilinear form
\begin{equation}
  a(u,v) = \int_\Omega \Delta u \Delta v \, dx,
\end{equation}
obtained by integrating by parts twice, presents analytic challenges (it fails to satisfy a G{\aa}rding-type inequality because  $a(v, v)$ vanishes for all harmonic functions).  

Following \citet{BreSco}, we employ the bilinear form
\begin{equation}
  a(u, v) = \int_{\Omega} \Delta u \Delta v - \left(1-\nu \right)
  \left( 2 u_{xx} v_{yy} + 2 u_{yy} v_{xx} - 4 u_{xy} v_{xy} \right)
  \dx,
  \label{eq:aplate}
\end{equation}
where $0 \leq \nu \leq \tfrac{1}{2}$ is the plate's Poisson ratio.
This form is actually coercive (modulo linear polynomials) for $0 < \nu < 1$.
The terms in $a$ multiplied by $(1-\nu)$ cancel out under integration by parts.
So, subject to clamped boundary conditions, the problem is well-posed.
Optimal error estimates can be shown in $H^1$ and $H^2$, although it is known
that obtaining optimal (third order) estimates in $L^2$ is quite challenging for spaces that contain only quadratic polynomials~\cite{brenner2005c,mu2019development}.

We used this example to validate our implementation of non-macro $C^1$ elements
in~\cite{finat-zany}, and we repeat this experiment for our
newly-implemented macroelements.  Macroelements
facilitate the strong enforcement of the clamped boundary conditions, in
contrast with~\cite{finat-zany}, where the supersmooth elements with second
derivative nodes required us to enhance the bilinear form with Nitsche-type
terms.
However, other choices of boundary conditions (for example $u=0=\tfrac{\partial^2u}{\partial n^2}$) might still require a Nitsche-type approach.

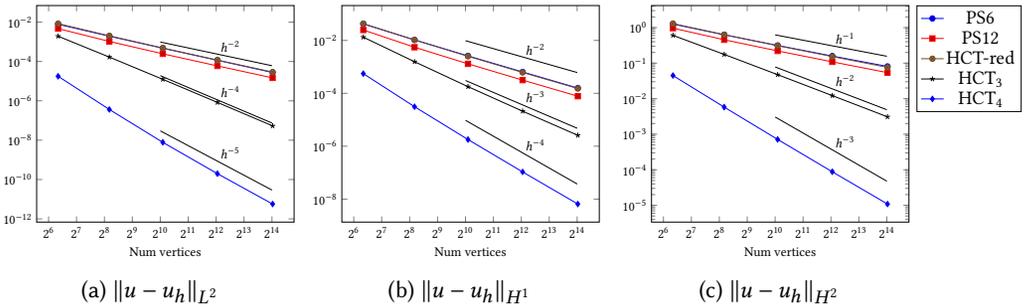
\begin{figure}[htbp]
 \begin{subfigure}[t]{0.29\textwidth}
   \begin{tikzpicture}[scale=0.5]
      \begin{loglogaxis}[xlabel={Num DOFs},
          legend style={font=\LARGE}, legend pos = outer north east,
          tick label style={font=\Large},
          label style={font=\Large}          
         ]
       \addplot table [x=NDOF, y=ErrorL2, col sep=comma]{code/biharmonic//data/biharmonic.PS6.2.csv};   %\addlegendentry{PS6}
       \addplot table [x=NDOF, y=ErrorL2, col sep=comma]{code/biharmonic//data/biharmonic.PS12.2.csv};   %\addlegendentry{PS12}
        \addplot table [x=NDOF, y=ErrorL2, col sep=comma]{code/biharmonic//data/biharmonic.HCT-red.3.csv};   %\addlegendentry{HCT-red}
        \addplot table [x=NDOF, y=ErrorL2, col sep=comma]{code/biharmonic//data/biharmonic.HCT.3.csv};       %\addlegendentry{$\mathrm{HCT}_3$}
        \addplot table [x=NDOF, y=ErrorL2, col sep=comma]{code/biharmonic//data/biharmonic.HCT.4.csv};       %\addlegendentry{$\mathrm{HCT}_4$}     
        \addplot [domain=1E4:1E5] {3E-4/pow(x/1E4,2/2)} node[above, yshift=-2pt, midway, anchor=south west] {$h^{2}$};
        %% \addplot [domain=1E4:1E5] {100/pow(x,3/2)} node[above, yshift=-1pt, midway, anchor=south west] {$h^{-3}$};
        \addplot [domain=1E4:1E5] {1E-5/pow(x/1E4,4/2)} node[above, yshift=-2pt, midway, anchor=south west] {$h^{4}$};
        \addplot [domain=1E4:1E5] {1E-7/pow(x/1E4,5/2)} node[above, yshift=-2pt,  midway, anchor=south west] {$h^{5}$};
   \end{loglogaxis}
  \end{tikzpicture}
  \caption{$\|u - u_h\|_{L^2}$}
 \end{subfigure}
 \begin{subfigure}[t]{0.29\textwidth}
   \begin{tikzpicture}[scale=0.5]
     \begin{loglogaxis}[xlabel={Num DOFs},
         legend style={font=\LARGE}, legend pos = outer north east,
          tick label style={font=\Large},
          label style={font=\Large}         
         ]
       \addplot table [x=NDOF, y=ErrorH1, col sep=comma]{code/biharmonic//data/biharmonic.PS6.2.csv};   %\addlegendentry{PS6}
       \addplot table [x=NDOF, y=ErrorH1, col sep=comma]{code/biharmonic//data/biharmonic.PS12.2.csv};   %\addlegendentry{PS12}
        \addplot table [x=NDOF, y=ErrorH1, col sep=comma]{code/biharmonic//data/biharmonic.HCT-red.3.csv};   %\addlegendentry{HCT-red}
        \addplot table [x=NDOF, y=ErrorH1, col sep=comma]{code/biharmonic//data/biharmonic.HCT.3.csv};       %\addlegendentry{$\mathrm{HCT}_3$}
        \addplot table [x=NDOF, y=ErrorH1, col
          sep=comma]{code/biharmonic//data/biharmonic.HCT.4.csv};
        \addplot [domain=1E4:1E5] {1.5E-3/pow(x/1E4,2/2)} node[above, yshift=-2pt, midway, anchor=south west] {$h^{2}$};
        \addplot [domain=1E4:1E5] {1.5E-4/pow(x/1E4,3/2)} node[above, yshift=-2pt, midway, anchor=south west] {$h^{3}$};
        \addplot [domain=1E4:1E5] {1E-5/pow(x/1E4,4/2)} node[above, yshift=-2pt,  midway, anchor=south west] {$h^{4}$};
   \end{loglogaxis}
  \end{tikzpicture}
  \caption{$\|u - u_h\|_{H^1}$}
 \end{subfigure} 
 \begin{subfigure}[t]{0.29\textwidth}
   \begin{tikzpicture}[scale=0.5]
     \begin{loglogaxis}[xlabel={Num DOFs},
         legend style={font=\LARGE}, legend pos = outer north east,
          tick label style={font=\Large},
          label style={font=\Large}         
         ]
       \addplot table [x=NDOF, y=ErrorH2, col sep=comma]{code/biharmonic//data/biharmonic.PS6.2.csv};   \addlegendentry{PS6}
       \addplot table [x=NDOF, y=ErrorH2, col sep=comma]{code/biharmonic//data/biharmonic.PS12.2.csv};   \addlegendentry{PS12}       
        \addplot table [x=NDOF, y=ErrorH2, col sep=comma]{code/biharmonic//data/biharmonic.HCT-red.3.csv};   \addlegendentry{HCT-red}
        \addplot table [x=NDOF, y=ErrorH2, col sep=comma]{code/biharmonic//data/biharmonic.HCT.3.csv};       \addlegendentry{$\mathrm{HCT}_3$}
        \addplot table [x=NDOF, y=ErrorH2, col
          sep=comma]{code/biharmonic//data/biharmonic.HCT.4.csv};
        \addlegendentry{$\mathrm{HCT}_4$}
        \addplot [domain=1E4:1E5] {3E-1/pow(x/1E4,1/2)} node[above, yshift=-2pt, midway, anchor=south west] {$h^{1}$};
        \addplot [domain=1E4:1E5] {5E-2/pow(x/1E4,2/2)} node[above, yshift=-2pt, midway, anchor=south west] {$h^{2}$};
        \addplot [domain=1E4:1E5] {2E-3/pow(x/1E4,3/2)} node[above, yshift=-2pt, midway, anchor=south west] {$h^{3}$};
        %% \addplot [domain=2^10:2^14] {10/pow(x,4/2)} node[above, yshift=-1pt,  midway, anchor=south west] {$h^{-4}$};
        %% \addplot [domain=2^6:2^10] {10/pow(x,5/2)} node[below, yshift=-2pt, midway] {$h^{-5}$};
   \end{loglogaxis}
  \end{tikzpicture}
  \caption{$\|u - u_h\|_{H^2}$}
 \end{subfigure}
 \caption{Error in solving biharmonic equation on a perturbed $N \times N$ mesh.}
  \Description{Convergence plots for biharmonic discretizations. The error in $L^2$, $H^1$, and $H^2$ norms are plotted against the number of degrees of freedom.} 
 \label{fig:biharmconv}  
\end{figure}

\Cref{fig:biharmconv} plots the error versus mesh refinement, where we take a
coarse mesh of the unit square, slightly perturb the internal vertices, and
then take uniform refinements.  Here, we see that HCT and its higher-order
variants provide more accurate results than the lower-order.
We note that the Powell-Sabin and reduced HCT elements give second order $H^1$ and first-order $H^2$ convergence rates, but only second order in $L^2$,
in accordance with the cited literature.

\begin{figure}[htbp]
    \begin{tikzpicture}
      \pgfplotstableread[col sep=comma,]{code/biharmonic/data/biharmonicflops.csv}\datatable
      \begin{axis}[ybar stacked,
          ylabel near ticks,
          %ymode = log,
        legend pos=north west,
        %% ymin=-1e4, ymax=0.5e4,
        legend cell align=left,
        legend style={draw=none},
        xtick=data,
        xticklabels from table={\datatable}{Name}, ylabel={FLOPs},
        xticklabel style={xshift=-10pt, rotate=60}]
        \addplot+ table [x expr=\coordindex, y={Cell Flops}]{\datatable};
        \legend{FLOPs}
      \end{axis}
    \end{tikzpicture}
    \caption{FLOP count for evaluating the element-level kernel for the biharmonic operator.  This count assumes reference basis elements are pre-tabulated and then includes the cost of transforming basis functions and their derivatives and performing integration.}
    \Description{Bar chart showing the number of floating point operations required to form the element matrices for various biharmonic elements.}     
    \label{fig:biharmflops}
\end{figure}
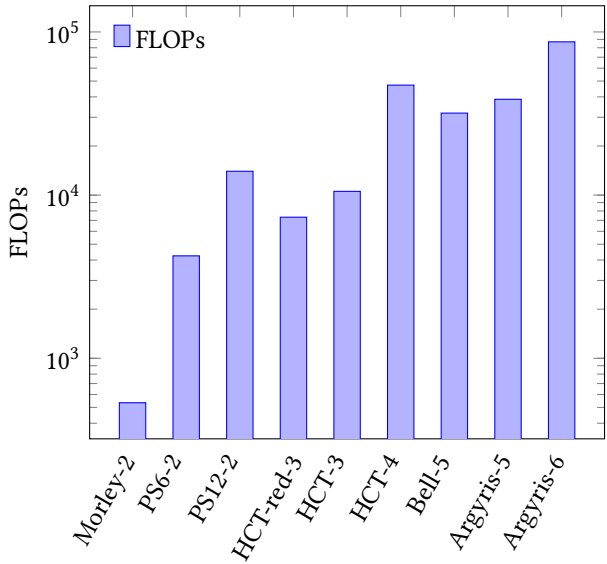

The macroelements give lower orders of accuracy, but have fewer global degrees
of freedom and lower polynomial degree than the Bell and Argyris elements.
However, because they use piecewise polynomials, we must integrate over each
subcell, so it is also interesting to compare the work required to form each
local stiffness matrix.  \Cref{fig:biharmflops} shows the \lstinline{tsfc}-reported FLOP counts for forming the local matrix over a single element, comparing our macroelements to
some
classical non-macroelements (the quadratic Morley element, plus Bell and two
degrees of Argyris).  The Powell-Sabin and cubic HCT elements indeed give a
middle ground between the inexpensive but low-order Morley element and the
higher-order polynomial $C^1$ elements.  We note that the $\HCT_4$ element is
actually more expensive than the Bell and quintic Argyris elements.  Although
the quintic Argyris element is more accurate and slightly lower-cost than the
$\HCT_4$ element, we note that some applications may benefit from the lack of
higher derivatives at the vertices.

\begin{figure}[htbp]
  \phantom{.}\hfill
  \begin{subfigure}[t]{0.15\textwidth}
    \centering
    \includegraphics[height=0.95\textwidth]{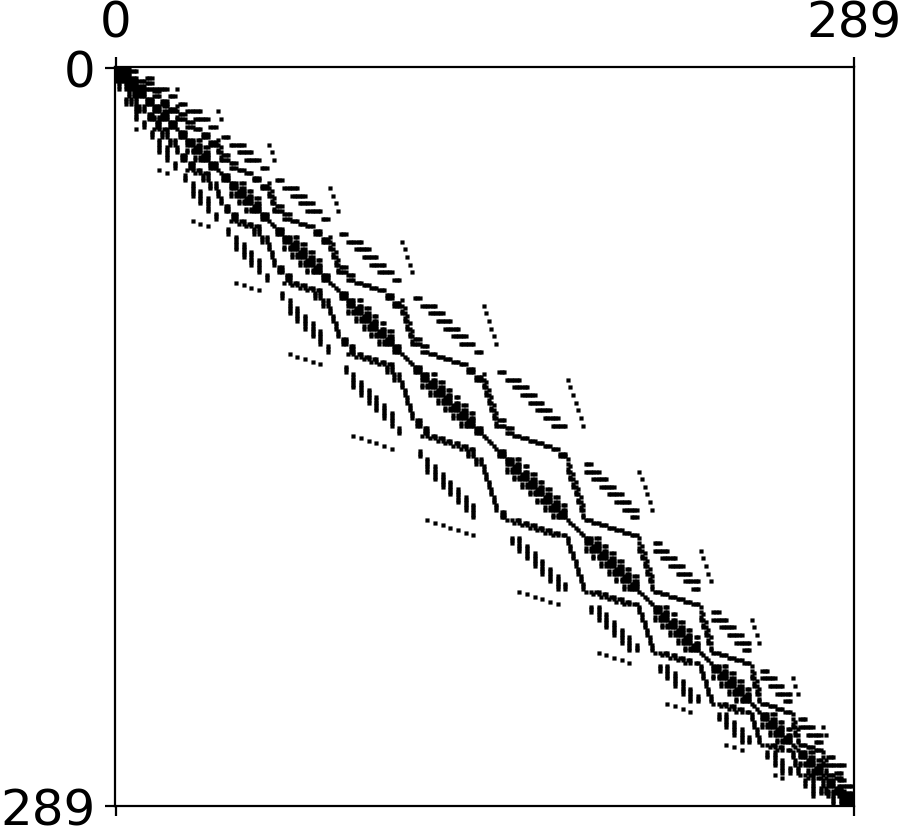}
    \caption{Morley}
  \end{subfigure}
  \hfill
  \begin{subfigure}[t]{0.15\textwidth}
    \centering
    \includegraphics[height=0.95\textwidth]{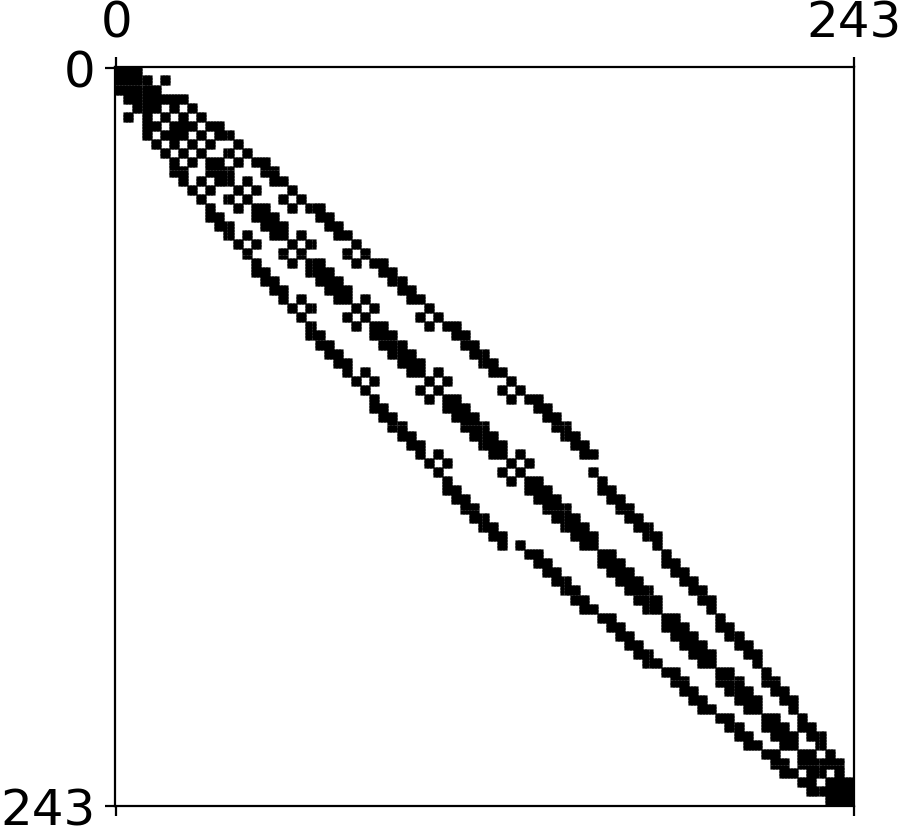}
    \caption{PS6}
  \end{subfigure}
  \hfill
  \begin{subfigure}[t]{0.15\textwidth}  
    \includegraphics[height=0.95\textwidth]{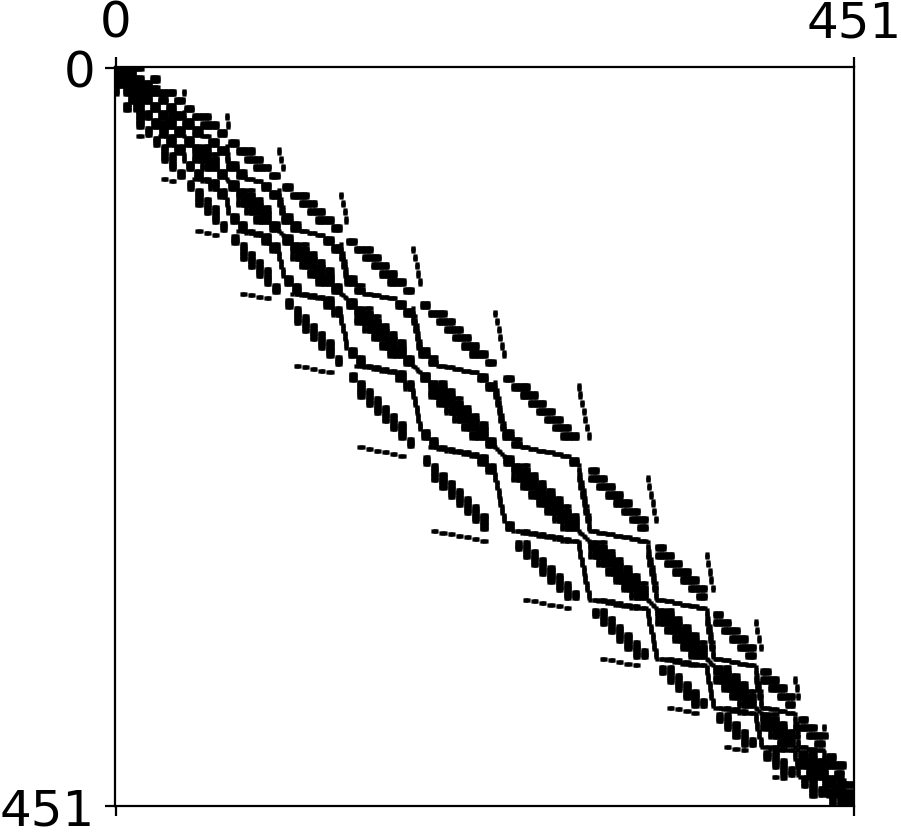}
    \caption{PS12}
  \end{subfigure}  
  \hfill
  \begin{subfigure}[t]{0.15\textwidth}
    \centering
    \includegraphics[height=0.95\textwidth]{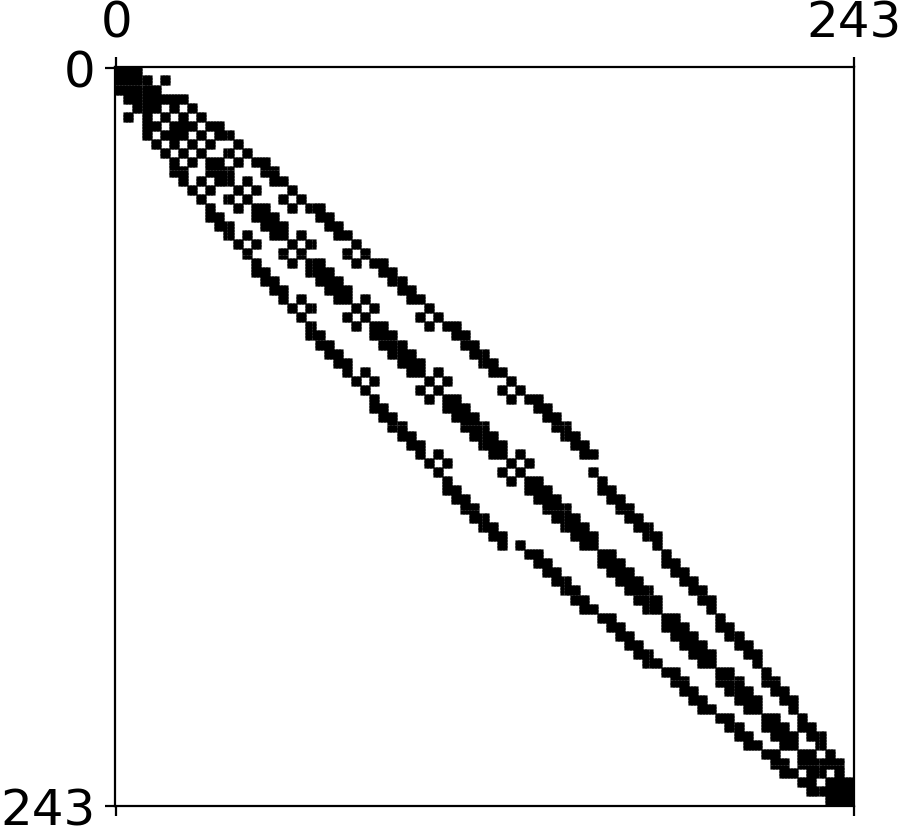}
    \caption{HCT-red}
  \end{subfigure}
  \hfill\phantom{.}

  \phantom{.}\hfill
  \begin{subfigure}[t]{0.15\textwidth}
    \centering
    \includegraphics[height=0.95\textwidth]{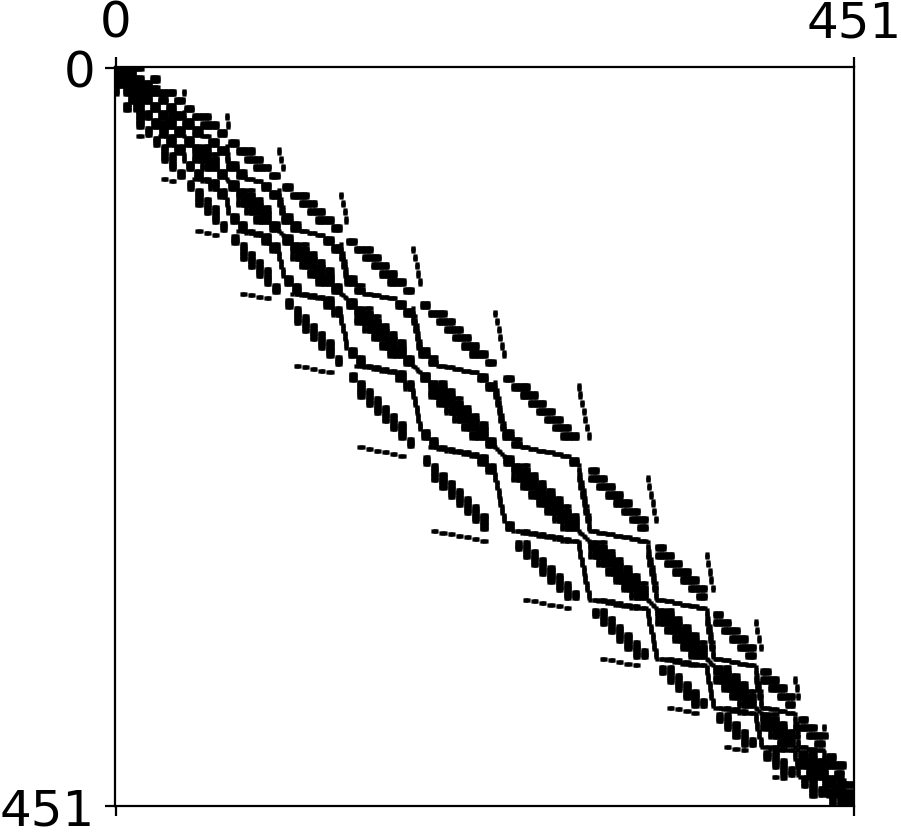}
    \caption{HCT-3}
  \end{subfigure}
  \hfill
  \begin{subfigure}[t]{0.15\textwidth}
    \centering
    \includegraphics[height=0.95\textwidth]{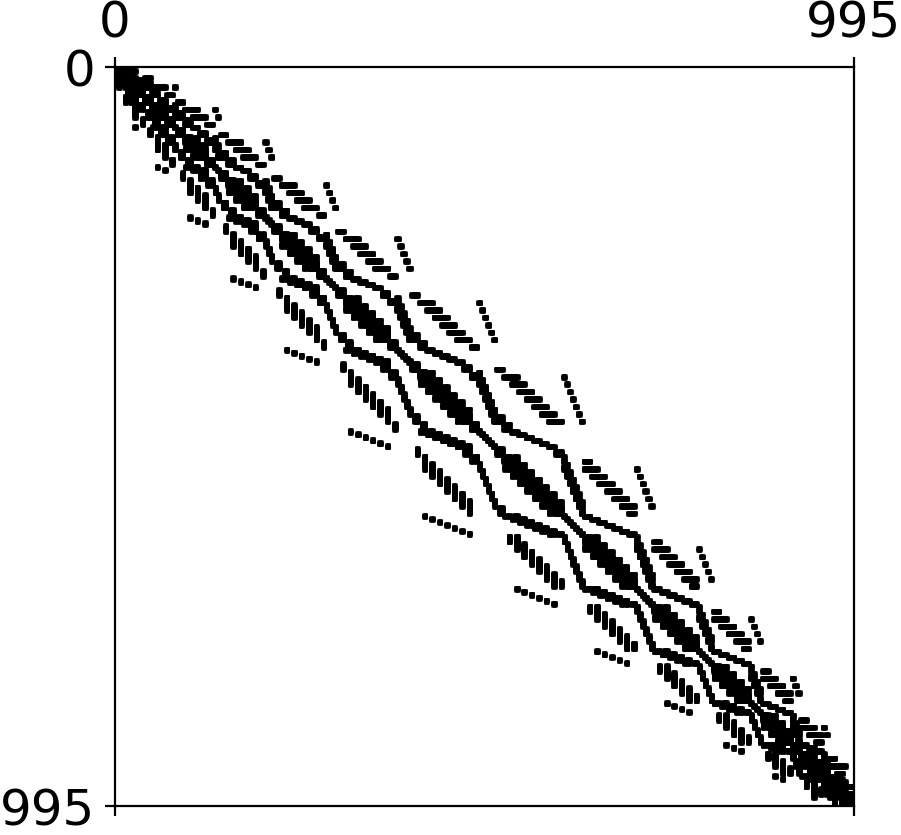}
    \caption{HCT-4}
  \end{subfigure}
  \hfill
  \begin{subfigure}[t]{0.15\textwidth}
    \centering
    \includegraphics[height=0.95\textwidth]{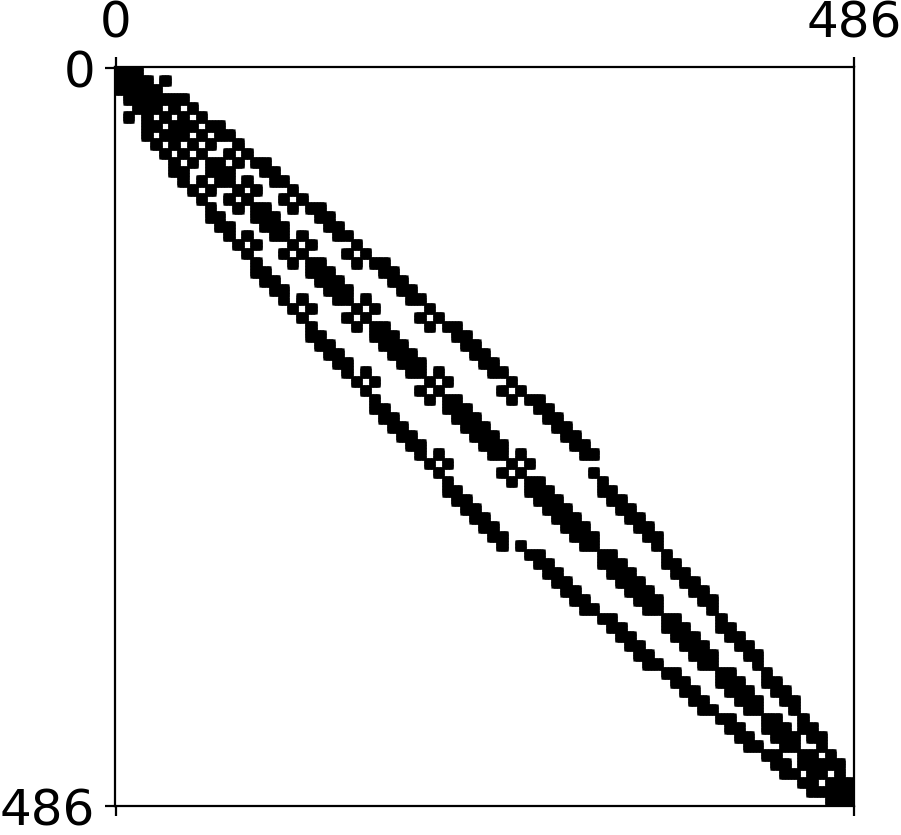}
    \caption{Bell}
  \end{subfigure}
  \hfill
  \begin{subfigure}[t]{0.15\textwidth}
    \centering
    \includegraphics[height=0.95\textwidth]{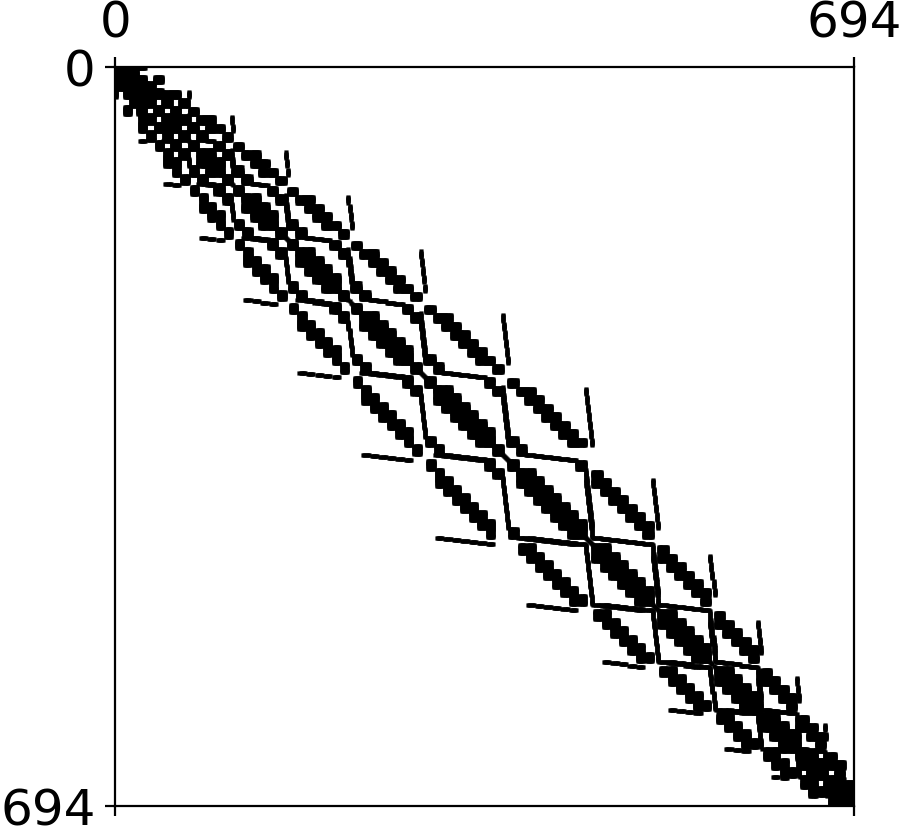}
    \caption{Argyris}
  \end{subfigure}    
  \hfill\phantom{.}
  \caption{Sparsity patterns of the discrete biharmonic operator on a regular $8
    \times 8$ mesh using various $H^2$ elements.}
  \label{fig:biharmspy}
  \Description{Plots of the sparsity patterns for the biharmonic operator on an $8 \times 8$ mesh using  various macroelements.}
\end{figure}

In addition the FLOP count for building a single element matrix, we report the number of nonzeros per row on an $8 \times 8$ mesh for several of our macroelements.  For comparison, we include the classical non-macro Morley, Bell, and Argyris elements.  The Morley element has only one degree of freedom for each edge and vertex, so leads to the sparsest matrix.  In contrast, the Bell and Argyris elements have six degrees of freedom per vertex (up to second derivatives) and lead to the densest matrices.  The quadratic Powell-Sabin-6 spline with only vertex degrees of freedom actually has a smaller but denser matrix than Morley.
Figure~\ref{fig:biharmspy} shows the sparsity pattern for several of our biharmonic elements on an $8 \times 8$ mesh.  We note that, because they have the same global degrees of freedom, the PS6 and HCT-red elements lead to the same sparsity patterns, as do PS12 and HCT3.  Moreover, PS6 and HCT-red as well as the Bell element have only vertex degrees of freedom.  Therefore, they have similar sparsity patterns, differing only in sizes of each block.

\begin{table}
  \begin{tabular}{ccc}
    Element & Rows & Nonzeros per row \\ \hline
    Morley & 289 & 10.63 \\
    PS6 & 243 & 18.41 \\
    PS12 & 451 & 22.73 \\
    HCT-red & 243 & 18.41 \\
    HCT3 & 451 & 22.73 \\
    HCT4 & 995 & 32.56 \\
    Bell & 486 & 36.81 \\
    Argyris & 694 & 41.02
  \end{tabular}
  \caption{The biharmonic equation was discretized on an $8 \times 8$ mesh, divided into right triangles.  This table shows the total number of degrees of freedom (rows in the matrix) and the average sparsity over all rows of the matrix.}
  \label{table:biharmonicsparsity}
\end{table}

\begin{figure}[htbp]
  \begin{subfigure}[t]{0.45\textwidth} \centering
    \begin{tikzpicture}[scale=0.6]
      \begin{loglogaxis}[xlabel=Num DOFs, legend style={font=\LARGE},
          tick label style={font=\Large},
          label style={font=\Large}, ylabel={Time (s)}
         ]
        \addplot table[x=NDOF, y=SNESSolve, col sep=comma]{code/biharmonic/data/biharmonic.PS6.2.csv};
        \addplot table[x=NDOF, y=SNESSolve, col sep=comma]{code/biharmonic/data/biharmonic.PS12.2.csv};
        \addplot table[x=NDOF, y=SNESSolve, col sep=comma]{code/biharmonic/data/biharmonic.HCT-red.3.csv};
        \addplot table[x=NDOF, y=SNESSolve, col sep=comma]{code/biharmonic/data/biharmonic.HCT.3.csv};
        \addplot table[x=NDOF, y=SNESSolve, col sep=comma]{code/biharmonic/data/biharmonic.HCT.4.csv};
        \addplot [domain=4E3:1E5] {1.75E-2 * pow(x/4E3,1)} node[above, xshift=16pt, yshift=-8pt, midway] {$\mathcal{O}(N)$};        
      \end{loglogaxis}
    \end{tikzpicture}
    \caption{Total run-time}
    \label{fig:biharmonicruntime}
  \end{subfigure}
  \begin{subfigure}[t]{0.45\textwidth} \centering
    \begin{tikzpicture}[scale=0.6]
      \begin{semilogxaxis}[xlabel=Num DOFs, legend
          style={font=\LARGE},legend pos = outer north east,
          tick label style={font=\Large},
          label style={font=\Large}, ymin=0, ymax=1,
          ylabel={Assembly fraction}
        ]
        \addplot table[x=NDOF, y expr= \thisrow{SNESJacobianEval} / \thisrow{SNESSolve}, col sep=comma]{code/biharmonic/data/biharmonic.PS6.2.csv}; \addlegendentry{PS6};
        \addplot table[x=NDOF, y expr= \thisrow{SNESJacobianEval} / \thisrow{SNESSolve}, col sep=comma]{code/biharmonic/data/biharmonic.PS12.2.csv}; \addlegendentry{PS12};
        \addplot table[x=NDOF, y expr= \thisrow{SNESJacobianEval} / \thisrow{SNESSolve}, col sep=comma]{code/biharmonic/data/biharmonic.HCT-red.3.csv}; \addlegendentry{HCT-red};
        \addplot table[x=NDOF, y expr= \thisrow{SNESJacobianEval} / \thisrow{SNESSolve}, col sep=comma]{code/biharmonic/data/biharmonic.HCT.3.csv}; \addlegendentry{HCT3};
        \addplot table[x=NDOF, y expr= \thisrow{SNESJacobianEval} / \thisrow{SNESSolve}, col sep=comma]{code/biharmonic/data/biharmonic.HCT.4.csv}; \addlegendentry{HCT4};
    \end{semilogxaxis}
    \end{tikzpicture}
    \caption{Fraction of time in assembly}
    \label{fig:biharmonicassemblyfraction}    
  \end{subfigure}
  \caption{Timing results versus the number of degrees of freedom for the two-dimensional biharmonic equation with various elements.  We show the total solve time (to assemble, factor, and solve) each system on the left and the fraction of that time spent in assembly on the right.}
  \Description{Plots showing the run-time and ratio of time spent in assembly to sparse direct factorization in some two-dimensional biharmonic calculations.}
  \label{fig:timingbiharmonic}
\end{figure}

Figure~\ref{fig:timingbiharmonic} shows the performance of solving the system with the sparse LU factorization from MUMPS.  Figure~\ref{fig:biharmonicruntime} shows the total time, including assembly, factorization, and backsolve, as a function of the total number of degrees of freedom in the mesh.  LU factorization is quite effective on a single core in 2D, only slightly superlinear.  Figure~\ref{fig:biharmonicassemblyfraction} shows the fraction of the total time spent in matrix assembly.  On the coarser meshes, we spend about 20\% of the time in assembly, while this drops to about 10\% on the finer meshes where factorization is relatively more expensive.

\section{Conclusions and future work}
\label{sec:conc}
Here, we have developed a robust, general framework for the construction of reference bases for macroelements in FIAT and their integration into Firedrake.
This appears to be the first such general-purpose implementation, and we are able to evaluate a representative suite of classical and modern elements for incompressible flow and biharmonic problems.
At the same time, these results suggest the need for much ongoing work.
For one, inter-grid transfers to enable multigrid algorithms for macroelement spaces are not fully understood, with mathematical and practical questions remaining.  Also, optimizing element-level calculation for macroelements, whether by use of sparse arrays or iterating over subcells with indirection, remains an open issue.
Finally, Firedrake supports sum-factorization of matrix-vector products on tensor product domains~\cite{homolya2017exposing}, and combining this with 1D ISO-type macroelements for preconditioning as in~\cite{pazner2023low} should be possible in a quite general setting.

\appendix
\section{Reproducibility}\label{sec:reproducibility}
The exact version of Firedrake used, along with scripts employed for the generation of numerical data is archived on Zenodo \citep{zenodo-macro}. 

\bibliographystyle{ACM-Reference-Format}
\bibliography{references}

\end{document}